\documentclass[11pt,english,english,english]{amsart}
\usepackage[latin1]{inputenc}
\usepackage{geometry}
\usepackage{setspace}
\usepackage{amssymb}

\makeatletter

\providecommand{\LyX}{L\kern-.1667em\lower.25em\hbox{Y}\kern-.125emX\@}

 \theoremstyle{plain}    
 \newtheorem{thm}{Theorem}[section]
 \numberwithin{equation}{section} 
 \numberwithin{figure}{section} 
 \theoremstyle{remark}    
 \newtheorem*{acknowledgement*}{Acknowledgement} 
 \theoremstyle{plain}    
 \newtheorem{cor}[thm]{Corollary} 
 \theoremstyle{plain}    
 \newtheorem{lem}[thm]{Lemma} 
 \theoremstyle{plain}    
 \newtheorem{prop}[thm]{Proposition} 
 \theoremstyle{plain}    
 \newtheorem*{thm*}{Theorem} 

\DeclareMathOperator{\aut}{Aut}
\DeclareMathOperator{\homeo}{Homeo}

\DeclareMathOperator{\diam}{diam}

\newcommand{\trans}{\textbf{T}}
\newcommand{\ttrans}{\textbf{T}^*}

\newcommand{\wm}{\textbf{WM}}
\newcommand{\sm}{\textbf{SM}}

\newcommand{\cgr}{\all_\textbf{Gr}}

\newcommand{\all}{\mathcal{S}}
\newcommand{\mn}{\textbf{Min}}
\newcommand{\ue}{\textbf{U}}

\newcommand{\cantor}{K}
\newcommand{\Cantor}{\mathcal{K}}

\usepackage{babel}
\makeatother

\usepackage{babel}
\makeatother

\usepackage{babel}
\makeatother
\begin{document}

\title{Genericity in topological dynamics}

\author{Michael Hochman}


\curraddr{Einstein Institute of Mathematics, Edmond J. Safra Campus, Givat
Ram, The Hebrew University of Jerusalem, Jerusalem 91904, Israel }

\email{mhochman@math.huji.ac.il}

\begin{abstract}
We study genericity of dynamical properties in the space of homeomorphisms
of the Cantor set and in the space of subshifts of a suitably large
shift space. These rather different settings are related by a Glasner-King
type correspondence: genericity in one is equivalent to genericity
in the other. 

By applying symbolic techniques in the shift-space model we derive
new results about genericity of dynamical properties for transitive
and totally transitive homeomorphisms of the Cantor set. We show that
the isomorphism class of the universal odometer is generic in the
space of transitive systems. On the other hand, the space of totally
transitive systems displays much more varied dynamics. In particular,
we show that in this space the isomorphism class of every Cantor system
without periodic points is dense, and the following properties are
generic: minimality, zero entropy, disjointness from a fixed totally
transitive system, weak mixing, strong mixing, and minimal self joinings.
The last two stand in striking contrast to the situation in the measure-preserving
category. We also prove a correspondence between genericity of dynamical
properties in the measure-preserving category and genericity of systems
supporting an invariant measure with the same property.
\end{abstract}
\maketitle
\pagestyle{myheadings}

\markboth{Michael Hochman}{Genericity in topological dynamics}

\section{\label{sec:Introduction}Introduction}

\subsection{\label{sub:Genericity}Genericity}

A question of some interest is, what does a {}``typical'' dynamical
system look like? To make this question precise one first fixes a
parameterization of dynamical systems with an appropriate complete
metric, and says that a property is \emph{generic} if the set of systems
satisfying the property is residual (contain a dense $G_{\delta }$),
and is \emph{exotic} if it is meagre (its complement is residual).

This question goes back to the early twentieth century, at least as
far as the work of Oxtoby and Ulam \cite{OU41}, who showed that in
dimension $2$ and up, a generic volume-preserving homeomorphism of
a cube (and other {}``nice'' manifolds) is ergodic. A few years
later, Halmos \cite{H44,Hal-mixing,Hal-ergodicth} showed that in
the space of automorphisms of a Lebesgue space, a generic automorphism
is ergodic and weak mixing. In contrast, Rohlin \cite{Rk-mixing}
showed that the strongly mixing automorphisms are exotic (this established
the existence of weakly but not strongly mixing systems before any
explicit examples were available). Since then the status of many other
properties in this category have been established. A parallel program
has been carried out for volume preserving homeomorphisms of manifolds,
and a certain unification has been achieved between these two categories
\cite{CP,AP02}. Similar questions have been studied in the smooth
category, where there are many open questions.

In the present paper we will be concerned with the category of topological
dynamics, where there has recently been renewed interest in questions
of genericity \cite{GW-TRproperty,BDK,AGW06,AHK03}. By a \emph{topological
dynamical system} we mean a pair $(X,\varphi )$ where $X$ is a compact
metric space and $\varphi :X\rightarrow X$ is a homeomorphism. Our
results concern two parameterizations of dynamical systems which are,
somewhat surprisingly, closely related: the space of subsystems of
a suitably large shift space, which is universal, and the group of
homeomorphisms of the Cantor set. We will mainly work in the first
of these, where symbolic techniques can be applied, but we begin the
discussion with the second, which is the better known of the two,
and where some interesting new phenomena have recently come to light.

\subsection{\label{sub:homeo-cantor}New Results for Homeomorphisms of the Cantor
Set}

Much work on genericity has focused on \emph{Cantor systems}, that
is, systems whose phase space is the Cantor set $\cantor $. Their
classical parametrization is the space of homeomorphisms of $\cantor $,
denoted \[
\mathcal{H}=\homeo (\cantor )\]
and topologized by the complete metric\[
d(\varphi ,\psi )=\max _{x\in \cantor }\left(d(\varphi (x),\psi (x))+d(\varphi ^{-1}(x),\psi ^{-1}(x))\right)\]
where on the right hand side $d(\cdot ,\cdot )$ is some fixed metric
for $\cantor $ (the same topology is induced by the usual metric
of uniform convergence, but note that that metric is not complete).
We identify each $\varphi \in \mathcal{H}$ with the Cantor dynamical
system $(\cantor ,\varphi )$. 

A striking and quite unexpected result about $\mathcal{H}$ was obtained
recently by Kechris and Rosendal \cite{KR04}, who showed that there
is a single Cantor system whose isomorphism class is a residual set
in $\mathcal{H}$. Hence, generically, there is only one Cantor system.
While Kechris and Rosendal's argument was nonconstructive, Akin, Glasner
and Weiss \cite{AGW06} describe this system explicitly, and it turns
out the its dynamics are quite degenerate and completely understood.
In particular, any question about genericity of a property in $\mathcal{H}$
is reduced to the question of whether this generic system satisfies
the property. 

For this reason we restrict attention to more dynamically interesting
subspaces of $\mathcal{H}$. One candidate is the space of transitive
systems. It turns out that here the periodic phenomena are dominant,
and again there is a single typical system:

\begin{thm}
\label{thm:transitive-systems}The space of transitive systems in
$\homeo (\cantor )$ is Polish, and the isomorphism class of the universal
odometer is generic there.
\end{thm}
The universal odometer is the unique Cantor system which factors onto
every finite cycle, and whose points are separated by these factors.
See section \ref{sub:Odometers} for details.

Much more interesting is the situation in the space of totally transitive
systems, which may be viewed as analogs of the aperiodic ergodic systems.
The absence of global periodicity leaves room for much more complex
dynamics:

\begin{thm}
\label{thm:totally-transitive-systems}The space of totally transitive
systems in $\homeo (\cantor )$ is Polish. Within this space, the
following hold: 
\begin{enumerate}
\item \label{thm:totally-transitive-systems-rokhlin}For a fixed Cantor
system without periodic points, the systems isomorphic to it are dense.
\item \label{thm:totally-transitive-systems-entropy}The zero-entropy, minimal
and uniquely ergodic systems are generic.
\item \label{thm:totally-transitive-systems-disjointness}For a fixed totally
transitive system, the systems disjoint from it are generic.
\item \label{thm:totally-transitive-systems-weak-mixing}The weakly mixing
systems are a dense $G_{\delta }$. 
\item \label{thm:totally-transitive-systems-strong-mixing}Strong mixing
is generic.
\item \label{thm:totally-transitive-systems-msj}Minimal self-joinings (and
hence semi-simplicity and primeness) are generic. 
\end{enumerate}
\end{thm}
It is interesting to compare these results with the situation in the
measure-preserving category, e.g. with the space of automorphisms
of a Lebesgue space. There are several points of similarity. (\ref{thm:totally-transitive-systems-rokhlin})
is an analogue of Halmos's classical theorem \cite{Hal-ergodicth}
that the isomorphism class of any aperiodic ergodic system is dense.
The statement about entropy in (\ref{thm:totally-transitive-systems-entropy})
is analogous to the situation in the measure preserving category.
Minimality and unique ergodicity may be viewed as {}``irreducibility''
conditions, and are somewhat analogous to ergodicity, which is likewise
generic. The measure-theoretic analogue of (\ref{thm:totally-transitive-systems-weak-mixing})
is due to Halmos \cite{Hal-ergodicth}. (\ref{thm:totally-transitive-systems-disjointness})
is an analog of a theorem of del Junco's \cite{dJ}, stating that
generically an automorphism is disjoint from a fixed ergodic system.
Note that in (\ref{thm:totally-transitive-systems-disjointness})
the fixed system need not be a Cantor system but can be any totally
transitive system.

On the other hand there are some striking differences. (\ref{thm:totally-transitive-systems-strong-mixing})
stands in contrast to the classical fact, due to Rohlin \cite{Rk-mixing},
that measure-theoretic strong mixing is exotic. (\ref{thm:totally-transitive-systems-msj})
should be compared to the fact, noted by del Junco \cite{dJ}, that
the measure preserving systems with minimal self joinings are exotic.
In particular, we note that primeness (i.e. having only trivial factors)
is generic in our setting but was recently proved to be exotic in
the measure-preserving case (Ageev \cite{A00}).

It should be noted that theorem \ref{thm:totally-transitive-systems}
isn't a result of the existence of a single generic system, as in
the space of transitive systems. This possibility is ruled out by
part (\ref{thm:totally-transitive-systems-disjointness}) of the theorem,
since a generic system would have to be disjoint from itself and hence
consist of one point, which is impossible for a Cantor system.

\subsection{\label{sub:shift-spaces-and-correspondence}The Shift-Space Model
and Correspondence Theorems}

Our motivation for this work was the question of genericity of general
topological systems, and in this regard $\mathcal{H}$ is not a particularly
good setting to investigate, since it represents only Cantor systems.
In general, the typical dynamics in spaces $\homeo (X)$ can depend
on $X$ in surprising and nontrivial ways. An example of this, due
to Glasner and Weiss \cite{GW-TRproperty}, is that zero entropy is
generic for homeomorphisms of the Cantor set, but infinite entropy
is generic for homeomorphisms of the Hilbert cube. These differences
are one reason it is desirable to have a universal model for dynamical
systems. 

Let $Q=[0,1]^{\aleph _{0}}$ denote the Hilbert cube and consider
the space $Q^{\mathbb{Z}}$ of bi-infinite sequences $x=(\ldots ,x(-1),x(0),x(1),\ldots )$
over $Q$, which is compact and metrizable in the product topology.
Let $\sigma $ be the shift homeomorphism of $Q^{\mathbb{Z}}$, i.e.
$(\sigma x)(i)=x(i+1)$. Let \[
\all =\{X\subseteq Q^{\mathbb{Z}}\, :\, X\textrm{ is closed, nonempty and }\sigma \textrm{-invariant}\}\]
be the space of subsystems of $Q^{\mathbb{Z}}$. In the Hausdorff
metric $\all $ is compact. We associate each $X\in \all $ with the
dynamical system $(X,\sigma |_{X})$, making $\all $ into a parametrization
of dynamical systems.

This parametrization is universal, that is, $\all $ contains members
of every isomorphism class of dynamical systems (see section \ref{sub:orbit-iamge-map}
below). However, it turns out that for purposes of studying genericity,
$\mathcal{S}$ and $\mathcal{H}$ are essentially the same. In order
to state this precisely we make the following definition: A dynamical
property $\mathbf{P}$ is a set of isomorphism classes of topological
dynamical systems. We abuse notation and write $(X,T)\in \mathbf{P}$
to indicate that the isomorphism class of $(X,T)$ is in $\mathbf{P}$.
If $\mathcal{U}$ is some parametrization of dynamical systems, we
write $\mathcal{U}_{\mathbf{P}}\subseteq \mathcal{U}$ for the set
of systems in $\mathcal{U}$ satisfying $\mathbf{P}$.

\begin{thm}
\label{thm:correspondence}(Correspondence theorem) Let $\textbf P$
be a dynamical property. Suppose that\renewcommand{\theenumi}{\alph{enumi}}
\begin{enumerate}
\item $\mathcal{H}_{\textbf P}$ is a $G_{\delta }$ in $\mathcal{H}$,
\item $\all _{\textbf P}$ is a $G_{\delta }$ in $\all $,
\item $\all _{\textbf P}$ contains a dense set of Cantor systems.
\end{enumerate}
\renewcommand{\theenumi}{\arabic{enumi}}Let $\textbf Q\subseteq \textbf P$
be a dynamical property. Then
\begin{enumerate}
\item $\textbf Q$ is generic in $\mathcal{H}_{\textbf P}$ if and only
if $\textbf Q$ is generic in $\all _{\textbf P}$.
\item If $\mathbf{Q}$ contains only Cantor systems, then $\textbf Q$ is
dense in $\mathcal{H}_{\textbf P}$ if and only if $\textbf Q$ is
dense in $\all _{\textbf P}$.
\end{enumerate}
\end{thm}
Conditions (a) and (b) are necessary for the statement of (1), since
in order to speak of genericity in $\all _{\mathbf{P}},\mathcal{H}_{\mathbf{P}}$
these spaces must be Polish, and a classical theorem of Alexandroff
states that this is equivalent to their being $G_{\delta }$'s (see
section \ref{sub:Baires-theory}). Condition (c) cannot be omitted
either. For example, let $\mathbf{P}$ be the property of being connected,
and let $\mathbf{Q}$ be the property of consisting of one point.
Then $\mathcal{H}_{\mathbf{P}}=\mathcal{H}_{\mathbf{Q}}=\emptyset $
and so density and genericity are satisfied trivially, whereas $\all _{\mathbf{Q}}$
is not dense in $\all _{\mathbf{P}}$.

Theorem \ref{thm:correspondence} is the analog of a similar result
in ergodic theory, where there are again two classical models of measure-preserving
systems. One is the automorphism group $\mathcal{A}$ of a Lebesgue
space, given the so-called coarse topology, which is analogous to
$\mathcal{H}$ (the phase space is fixed and the map varies). The
other is the space $\mathcal{M}$ of shift-invariant Borel probability
measures on $Q^{\mathbb{Z}}$, with the weak-{*} topology, which is
analogous to $\all $ (we fix the map and vary the subsystem/measure).
See section \ref{sub:Measure-Preserving-Systems} for more details.
Both these spaces are universal, and it has been shown by Glasner
and King \cite{GK-01law} and also by Rudolph \cite{R} that a property
is generic in one if and only if it is generic in the other. It is
worth mentioning that in ergodic theory there is a third universal
parametrization of dynamical systems, namely the space of transformations
which are orbit equivalent to a fixed ergodic transformation. Rudolph
\cite{R} has examined this space and shown that genericity of dynamical
properties there is equivalent to genericity in $\mathcal{A}$ and
$\mathcal{M}$. It remains to be explored to what extend there is
a topological analog of this.

Besides its intrinsic interest, our main use of the correspondence
theorem is in the proof of theorems \ref{thm:transitive-systems}
and \ref{thm:totally-transitive-systems}. We will verify the hypotheses
of the correspondence theorem for the properties $\textbf P$ of being
a Cantor system, being transitive, and being totally transitive systems.
Using symbolic techniques we will then be able to prove results on
the genericity of various properties in $\mathcal{S}$, and transfer
them to $\mathcal{H}$. We also can go the other way, and deduce from
Kechris and Rosendal's results that in $\mathcal{S}$ there is a single
generic system. 

An interesting correspondence also exists between genericity in the
measure-preserving category and the topological one. For a dynamical
property $\mathbf{P}$ in the measure-preserving category, we write
$\mathcal{M}_{\mathbf{P}}\subseteq \mathcal{M}$ for the set of invariant
Borel probability measures $\mu $ on $Q^{\mathbb{Z}}$ such that
the measure preserving system $(Q^{\mathbb{Z}},\sigma ,\mu )$ is
represented in $\mathbf{P}$. We have the following correspondence
principle:

\begin{thm}
\label{thm:inter-category-correspondence}Let $\mathbf{P}$ be a dynamical
property in the measure-preserving category and suppose that $\mathcal{M}_{\mathbf{P}}$
is a dense $G_{\delta }$ in $\mathcal{M}$. Let $\widetilde{\mathbf{P}}$
be the property (in the topological category) of supporting an invariant
measure from $\mathbf{P}$. Then $\widetilde{\mathbf{P}}$ is generic
in the space of totally transitive systems in $\all $ and $\mathcal{H}$.
\end{thm}
This theorem is similar to results of Alpern \cite{A78}, who showed
that under some conditions a $G_{\delta }$ property which is generic
for automorphisms of a Lebesgue space is also generic for volume-preserving
homeomorphisms of a manifold. The hypothesis in theorem \ref{thm:inter-category-correspondence},
that $\mathcal{M}_{\mathbf{P}}$ be a $G_{\delta }$, is unfortunate
(note that such an hypothesis is present also in Alpern's work). It
seems possible that the correspondence holds more generally, but we
do not know how to prove it.

Finally, an application of Glasner and King's techniques \cite{GK-01law}
give a zero-one law for {}``well-behaved'' dynamical properties.
Recall that a subset $A$ of a Polish space is a Baire set if it can
be written as $A=U\bigtriangleup M$, where $U$ is open and $M$
is first category.

\begin{thm}
\label{thm:zero-one-law}Let $\textbf P$ be a topological dynamical
property such that $\mathcal{H}_{\mathbf{P}}$ has the property of
Baire. Then either $\mathbf{P}$ is generic in the space of totally
transitive systems in $\mathcal{H}$ and $\all $, or it is exotic
there.
\end{thm}
Although our use of $\mathcal{S}$ is mainly for the study of Cantor
systems, we note again that it is a universal model, and thus contains
many other interesting subspaces, e.g. the space of all connected
systems. There are many interesting questions here which we have not
touched on. We discuss some of them in section \ref{sec:open-questions}.

\subsection{\label{sub:Notes-on-Notation}Organization and Notation }

The rest of this paper is organized as follows. In the next section
we give the basic definitions and notation. Section \ref{sec:correspondence}
is fairly independent and is devoted to the proof of the correspondence
between $\mathcal{H}$ and $\all $ and the zero-one law. Sections
\ref{sec:symbolic-approximation}-\ref{sec:prime-and-doubly-minimal-systems}
contain the remaining results: Section \ref{sec:symbolic-approximation}
contains some further definitions and symbolic machinery, section
\ref{sec:transitivity-and-periodic-approximation} deals with transitive
systems and proves theorem \ref{thm:transitive-systems}, section
\ref{sec:total-transitivity} deals with totally transitive systems
and proves parts (\ref{thm:totally-transitive-systems-rokhlin}) and
(\ref{thm:totally-transitive-systems-entropy}) of theorem \ref{thm:totally-transitive-systems}
as well as theorem \ref{thm:inter-category-correspondence}. Section
\ref{sec:Disjointness} proves part (\ref{thm:totally-transitive-systems-disjointness})
of theorem \ref{thm:totally-transitive-systems}. Section \ref{sec:mixing}
discusses parts (\ref{thm:totally-transitive-systems-weak-mixing})
and (\ref{thm:totally-transitive-systems-strong-mixing}) of that
theorem, and part (\ref{thm:totally-transitive-systems-msj}) is proved
in section \ref{sec:prime-and-doubly-minimal-systems}. Finally, in
section \ref{sec:open-questions} we outline some extensions and open
questions.

For the readers convenience, we conclude the introduction with a summary
of our main notational and typographical conventions. Further definitions
and notation appear in sections \ref{sec:Definitions} and \ref{sec:symbolic-approximation}.

\begin{center}\begin{tabular}{c|l}
$K$&
The Cantor set.\\
\hline 
$Q$&
The Hilbert cube, $[0,1]^{\aleph _{0}}$.\\
\hline 
$X,Y,Z\ldots $&
The phase space of dynamical systems.\\
\hline 
$x,y,z\ldots $&
Points in $X,Y,Z,\ldots $ etc. See section \ref{sub:Approximation-by-symbolic}
regarding subscripts.\\
\hline 
$x(i)$&
For $x\in Q^{\mathbb{Z}}$, the $i$-th coordinate of $x$.\\
\hline 
$S,T,f,g,\varphi ,\psi \ldots $&
Functions (usually continuous)\\
\hline 
$\sigma $&
The shift map on $Q^{Z}$\\
\hline 
$[(X,T)]$&
Isomorphism class of $(X,T)$.\\
\hline 
$\mathbf{P},\mathbf{Q},\ldots $&
Dynamical properties, i.e. sets of isomorphism classes.\\
\hline 
$\mathcal{H}$ , $\mathcal{H}_{\mathbf{P}}$&
The space of homeomorphisms of the Cantor set {[}satisfying $\mathbf{P}${]}.\\
\hline 
$\mathcal{S}$ , $\mathcal{S}_{\mathbf{P}}$&
The space of shift-invariant subsystems of $Q^{\mathbb{Z}}$ {[}satisfying
$\mathbf{P}${]}.\\
\hline 
$\mathcal{A}$ , $\mathcal{A}_{\mathbf{P}}$&
The automorphism group of a Lebesgue space {[}satisfying $\mathbf{P}${]}.\\
\hline 
$\mathcal{M}$ , $\mathcal{M}_{\mathbf{P}}$&
The shift-invariant probability measures on $Q^{\mathbb{Z}}$ {[}satisfying
$\mathbf{P}${]}.\\
\hline 
$d,d_{\infty }$&
Metrics, vary with context. See sections \ref{sub:homeo-cantor-defs},
\ref{sub:The-Hausdorf-Metric}, \ref{sub:Symbolic-Systems}, \ref{sub:Projection-into-symbolic-systems}.\\
\hline 
$\pi _{0}$&
Projection onto the $0$-th coordinate of $Q^{\mathbb{Z}}$\\
\hline 
$f_{T}^{*}$&
Orbit image map of $f$. See section \ref{sub:orbit-iamge-map}.\\
\hline 
$\theta $&
Pseudo-orbit-to-orbit map. See section \ref{sub:Approximation-by-symbolic}.\\
\hline 
$\tau _{Y}$&
Symbolic retract {}``towards'' $Y$. See section \ref{sub:Projection-into-symbolic-systems}.\\
\end{tabular}\end{center}

\begin{acknowledgement*}
This paper is part of the author's Ph.D. thesis, conducted under the
guidance of Professor Benjamin Weiss, whom I would like to thank for
all his support and advice. 
\end{acknowledgement*}

\section{\label{sec:Definitions}Definitions}

In this section we present some definitions and notation. For convenience,
we repeat here some of the definitions given in the introduction.
We warn the reader in advance that the letter $d$ will simultaneously
denote metrics on many spaces; which metric depends on the context.

\subsection{\label{sub:Dynamical-systems}Dynamical systems}

A \emph{topological dynamical system} is a pair $(X,T)$ where $X$
is a nonempty compact metric space and $T$ a homeomorphism of $X$.
The space $X$ is called the \emph{phase space}. We often write only
$X$ or $T$ in place of $(X,T)$. We will sometimes restrict to the
case where $X$ does not have isolated points; in this case $X$ is
said to be \emph{perfect}. 

A system $(Y,S)$ is a \emph{factor} of $(X,T)$ if there is a continuous
onto map $\varphi :X\rightarrow Y$ such that $\varphi T=S\varphi $.
Such a $\varphi $ is called a \emph{homomorphism} or \emph{factor
map} from $X$ to $Y$. If $\varphi $ is a homeomorphism it is called
an \emph{isomorphism,} and $(X,T),(Y,S)$ are said to be isomorphic.

A subset $X_{0}\subseteq X$ of a dynamical system $(X,T)$ is invariant
if $TX_{0}=T^{-1}X_{0}=X_{0}$. A closed, nonempty invariant subset
$X_{0}$ of $X$ defines a dynamical system by restricting $T$ to
$X_{0}$ and is called a \emph{subsystem} of $(X,T)$.

\subsection{\label{sub:homeo-cantor-defs}Homeomorphisms of the Cantor Set}

We denote by $\cantor $ the Cantor set. For some fixed metric $d$
on $\cantor $ we define\[
d(\varphi ,\psi )=\max _{x\in \cantor }\left(d(\varphi (x),\psi (x))+d(\varphi ^{-1}(x),\psi ^{-1}(x))\right)\]
for $\varphi ,\psi $ homeomorphisms of $\cantor $. We denote by\[
\mathcal{H}=\homeo (\cantor )\]
the space of homeomorphisms with the metric defined above, which is
complete.

\subsection{\label{sub:The-Hausdorf-Metric}The Hausdorff Metric}

Let $X$ be a compact metric space with metric $d$; as usual $B_{r}(x)$
is the open ball of radius $r$ around $x$. The space of nonempty,
closed subsets of $X$ is denoted by $2^{X}$, and the \emph{Hausdorff
metric} is defined on $2^{X}$ by \begin{eqnarray*}
d(Y_{0},Y_{1}) & = & \inf \left\{ \varepsilon \left|\begin{array}{c}
 \textrm{for }i=0,1\textrm{ and }y\in Y_{i}\textrm{ there is}\\
 y'\in Y_{1-i}\textrm{ with }d(y,y')<\varepsilon \end{array}
\right.\right\} 
\end{eqnarray*}
 for nonempty, closed subsets $Y_{0},Y_{1}\subseteq X$. With the
Hausdorff metric, $2^{X}$ is compact. This topology on $2^{X}$ can
be characterized as follows: If $Y_{n},Y\subseteq X$ are closed sets,
then $Y_{n}\rightarrow Y$ if and only if for any sequence $(y_{n})$
with $y_{n}\in Y_{n}$ the set of accumulation points of $(y_{n})$
is contained in $Y$, and every point in $Y$ arises in this way.

Note that if $(X,T)$ is a dynamical system then $T$ induces a homeomorphism
$\widetilde{T}$ of $2^{X}$, and the subsystems are precisely the
fixed points of $\widetilde{T}$. In particular the space of subsystems
is closed in the Hausdorff metric.

\subsection{\label{sub:The-space-of-Subshifts}The space of Subshifts}

Let $Q$ denote the Hilbert cube (the product of $[0,1]$ with itself
countably many times) with a fixed metric $d$. Let $Q^{\mathbb{Z}}$
be the space of bi-infinite sequences $(\ldots ,x(-1),x(0),x(1),\ldots )$
with $x(i)\in Q$, along with the product topology, which is also
compact and metrizable. To be concrete, we define a compatible metric
$d$ on $Q^{\mathbb{Z}}$ by \[
d(y,y')=\inf \{\varepsilon >0\, :\, d(y(i),y'(i))<\varepsilon \textrm{ for all }|i|\leq 1/\varepsilon \}\]
 for $y,y'\in Q^{\mathbb{Z}}$. Thus $y,y'$ are close in $Q^{\mathbb{Z}}$
if their coordinates agree well on a large block of indices around
zero. Note that we have arranged things so that if $d(y,y')<\varepsilon $
then $d(y(0),y'(0))<\varepsilon $. 

The \emph{shift $\sigma $} on $Q^{\mathbb{Z}}$ is the homeomorphism
defined by\[
\sigma (x)(n)=x(n+1)\]
 The dynamical system $(Q^{\mathbb{Z}},\sigma )$ is called the \emph{full
shift} on $Q$. 

The space of closed subsystems of $(Q^{\mathbb{Z}},\sigma )$ is \[
\all =\{X\subseteq Q^{Z}\, :\, X\neq \emptyset \textrm{ is closed and }\sigma -\textrm{invariant}\}\]
 As noted above, $\mathcal{S}$ is closed and compact in the Hausdorff
metric. 

We denote by $\pi _{0}$ the projection $Q^{\mathbb{Z}}\rightarrow Q$
onto the $0$-th coordinate, i.e. $\pi _{0}(x)=x(0)$. For $X\in \all $
the set $\pi _{0}(X)$ is called the \emph{cross-section} of $X$.

\subsection{\label{sub:orbit-iamge-map}Embedding Dynamical Systems in Shift
Spaces}

Let $(X,T)$ be a dynamical system and $x\in X$. The \emph{(full)
orbit} of $x$ is the set $\{T^{k}x\}_{k\in \mathbb{Z}}=\{\ldots ,T^{-1}x,x,Tx,T^{2}x,\ldots \}$.
Let $f:X\rightarrow Q$ be a continuous map. $f$ defines an \emph{orbit
picture} of $(X,T)$ in $Q^{\mathbb{Z}}$ by assigning to $x\in X$
the sequence of images under $f$ of its orbit: \[
f_{T}^{*}(x)=(f(T^{k}x))_{k\in \mathbb{Z}}\]
 One verifies that $f_{T}^{*}(Tx)=\sigma f_{T}^{*}(x)$, so $f_{T}^{*}$
is a factor map from $(X,T)$ onto its image. If $f$ is an embedding
of $X$ in $Q$ then $(f_{T}^{*}(X),\sigma )$ is isomorphic to $(X,T)$
as a dynamical system via the isomorphism $f_{T}^{*}$.

Since $Q$ has the property that any compact metric space can be embedded
in it, the previous discussion implies that any dynamical system $(X,T)$
can be embedded as a subsystem of $(Q^{\mathbb{Z}},\sigma )$; thus
$\all $ contains representatives of every isomorphism class of dynamical
systems.

\subsection{\label{sub:Dynamical-Properties}Dynamical Properties}

A \emph{dynamical property} $\textbf P$ is a family of isomorphism
classes of dynamical systems. For a system $(X,T)$ we write $[(X,T)]$
for its isomorphism class, although we shall write $(X,T)\in \textbf P$
instead of $[(X,T)]\in \textbf P$ (or just $X\in \textbf P$ or $T\in \mathbf{P}$
when $X$ or $T$ are understood). When considering a space of dynamical
systems such as $\mathcal{H}$ or $\all $ we identify $\textbf P$
with the subset of the space made up of those systems which have the
property $\textbf P$, and introduce the notation\begin{eqnarray*}
\mathcal{H}_{\textbf P} & = & \{\varphi \in \mathcal{H}\, :\, (\cantor ,\varphi )\in \textbf P\}\\
\all _{\textbf P} & = & \{X\in \all \, :\, (X,\sigma |_{X})\in \textbf P\}
\end{eqnarray*}
When we attribute topological properties such as openness or denseness
to $\textbf P$, we are actually referring to the sets $\all _{\textbf P}$
or $\mathcal{H}_{\textbf P}$. In particular we say $\textbf P$ is
dense, generic or exotic in $\all ,\mathcal{H}$ if $\all _{\textbf P},\mathcal{H}_{\textbf P}$
are respectively dense, residual or meager.

\subsection{\label{sub:Baires-theory}Baire Category}

We briefly review of the basic facts from Baire's category theory
that we will use. A good reference for this is Oxtoby's book \cite{Ox}.

A topological space $X$ is \emph{Polish} if there is a complete separable
metric $d$ on $X$ which induces the topology. A countable intersection
of open sets is called a $G_{\delta }$ set. A set containing a dense
$G_{\delta }$ is \emph{residual}; the complement of a residual set
is \emph{meagre} or \emph{first category} (any set which isn't first
category is second category, but this doesn't imply residuality).
Baire's theorem states that in a Polish space the intersection of
countably many dense open sets is dense; thus the intersection of
countably many residual sets is residual. The family of residual sets
forms a $\sigma $-filter on $X$; in this sense, residual sets are
the topological analogue of sets of full measure.

Alexandroff's classical theorem characterizes the Polish subsets of
a Polish space:

\begin{thm}
\label{thm:Alexandroff}(Alexandroff \cite[Theorem 12.1 and 12.3]{Ox})
For a Polish space $X$, a set $A\subseteq X$ is Polish if and only
if it is a $G_{\delta }$. 
\end{thm}

\subsection{\label{sub:Measure-Preserving-Systems}Measure Preserving Systems}

A measure preserving system is a quadruple $(X,\mathcal{F},\mu ,T)$
where $(X,\mathcal{F},\mu )$ is a standard probability space, $T:X\rightarrow X$
is bi-measurable (i.e. $T$ and $T^{-1}$ are measurable) and $T$
preserves $\mu $, i.e. $\mu (T^{-1}(A))=\mu (A)$ for all $A\in \mathcal{F}$.
Often $X$ will be a topological space, and then we always assume
that $\mathcal{F}$ is the completion of the Borel $\sigma $-algebra
with respect to a regular measure $\mu $. Two measure-preserving
systems $(X,\mathcal{F},\mu ,T)$ and $(Y,\mathcal{F}',\nu ,S)$ are
isomorphic if there is a measure-preserving invertible map $\pi :X\rightarrow Y$
such that $\pi T=S\pi $. 

Let $\lambda $ denote Lebesgue measure on $[0,1]$ and $\mathcal{L}$
the Lebesgue sets. Let\[
\mathcal{A}=\aut ([0,1],\mathcal{L},\lambda )\]
denote the set of measure-preserving automorphisms of $([0,1],\mathcal{L},\lambda )$.
Each $\varphi \in \mathcal{A}$ may be identified with the measure-preserving
system $([0,1],\mathcal{L},\lambda ,\varphi )$. This space is universal
for non-atomic measure preserving systems. Define a topology on $\mathcal{A}$
by $\varphi _{n}\rightarrow \varphi $ if $\varphi _{n}(A)\rightarrow \varphi (A)$
for all $A\in \mathcal{L}$. This topology is Polish.

A Borel probability measure $\mu $ on $Q^{\mathbb{Z}}$ is invariant
under the shift if $\mu (\sigma ^{-1}(A))=\mu (A)$ for every Borel
set $A\subseteq Q^{Z}$. Denote the space of shift-invariant measures
by $\mathcal{M}$. We may associate to $\mu \in \mathcal{M}$ the
measure preserving system $(Q^{\mathbb{Z}},\mathcal{F},\sigma ,\mu )$.
Identifying measures with positive linear functionals on $C(Q^{\mathbb{Z}})$
and using the Riesz representation theorem, we may equip $\mathcal{M}$
with the weak-{*} topology, which turns it into a compact metrizable
space. This space is universal for measure preserving systems.

\section{\label{sec:correspondence}The correspondence theorem and zero-one
laws}

\subsection{\label{sub:Formulation-and-Proof}Formulation and Proof Outline}

Our correspondence theorem is modeled after Glasner and King's result
\cite{GK-01law}, which states that a dynamical property in the measure
preserving category is generic in the automorphism group $\mathcal{A}$
of a Lebesgue space if and only if it is generic in the space $\mathcal{M}$
of shift invariant measures on $Q^{\mathbb{Z}}$ (the topologies on
these spaces were defined in the previous section). In our setting
the exact analogy of this would be that a dynamical property $\textbf P$
of topological systems is generic in $\mathcal{H}$ if and only if
it is generic in $\all $. Although this is true it is of limited
interest, since as we have already mentioned, generically there is
only one system in $\mathcal{H}$ up to isomorphism, so genericity
of $\mathbf{P}$ reduces to the question of whether this generic system
satisfy $\mathbf{P}$ or not. The version of the correspondence described
below is a relative one, asserting the if $\mathbf{P}$ is a dynamical
property satisfying certain conditions and $\mathbf{Q}\subseteq \mathbf{P}$,
then $\mathbf{Q}$ is generic in $\mathcal{H}_{\mathbf{P}}$ if and
only if $\mathbf{Q}$ is generic in $\all _{\mathbf{P}}$. We will
later apply this to the case where $\mathbf{P}$ is the class of transitive
or totally transitive systems.

We recall for convenience the formulation of theorem \ref{thm:correspondence}:

\begin{thm*}
(Correspondence theorem) Let $\textbf P$ be a dynamical property.
Suppose\renewcommand{\theenumi}{\alph{enumi}}
\begin{enumerate}
\item $\mathcal{H}_{\textbf P}$ is a $G_{\delta }$ in $\mathcal{H}$,
\item $\all _{\textbf P}$ is a $G_{\delta }$ in $\all $,
\item $\all _{\textbf P}$ contains a dense set of Cantor systems.
\end{enumerate}
\renewcommand{\theenumi}{\arabic{enumi}}Let $\textbf Q\subseteq \textbf P$
be a dynamical property. Then
\begin{enumerate}
\item $\textbf Q$ is generic in $\mathcal{H}_{\textbf P}$ if and only
if $\textbf Q$ is generic in $\all _{\textbf P}$.
\item If $\mathbf{Q}$ contains only Cantor systems, then $\textbf Q$ is
dense in $\mathcal{H}_{\textbf P}$ if and only if $\textbf Q$ is
dense in $\all _{\textbf P}$.
\end{enumerate}
\end{thm*}
\begin{proof}
The proof proceeds as follows. We first construct a Polish space $\mathcal{E}$
and a map $\beta :\mathcal{H}_{\textbf P}\times \mathcal{E}\rightarrow \all _{\textbf P}$
such that
\begin{itemize}
\item For all $\varphi \in \mathcal{H}_{\mathbf{P}}$ and $f\in \mathcal{E}$,
the system $(\beta (\varphi ,f),\sigma )$ is isomorphic to the system
$(\cantor ,\varphi )$. 
\item $\beta $ is a topological embedding.
\item The image of $\beta $ is dense, and furthermore if $X\in \all _{\mathbf{P}}$
is a Cantor system then $X=\lim X_{n}$ for a sequence of $X_{n}$
which are isomorphic to $X$ and contained in the image of $\beta $.
\end{itemize}
Before undertaking the construction of $\mathcal{E}$ and $\beta $,
which is somewhat involved, let us show how their existence proves
the theorem. Denote $\beta $'s image by $\mathcal{B}\subseteq \all _{\mathbf{P}}$.
We first claim  that $\mathcal{B}$ is a dense $G_{\delta }$ in $\all _{\mathbf{P}}$.
Indeed, it is dense by assumption, and since $\mathcal{H}_{\textbf P}\times \mathcal{E}$
is Polish and $\beta $ a topological embedding, $\mathcal{B}$ is
Polish and by Alexandroff's theorem is a $G_{\delta }$ in $\all _{\textbf P}$. 

Denote $\mathcal{B}_{\mathbf{Q}}=\mathcal{B}\cap \all _{\mathbf{Q}}$
(this may not be a $G_{\delta }$). We claim that $\beta ^{-1}(\mathcal{B}_{\mathbf{Q}})=\mathcal{H}_{\mathbf{Q}}\times \mathcal{E}$.
Indeed, this is because\begin{eqnarray*}
\beta ^{-1}(\mathcal{B}_{\mathbf{Q}}) & = & \{(\varphi ,f)\in \mathcal{H}_{\mathbf{P}}\times \mathcal{E}\, :\, (\beta (\varphi ,f),\sigma )\in \mathbf{Q}\}\\
 & = & \{(\varphi ,f)\in \mathcal{H}_{\mathbf{P}}\times \mathcal{E}\, :\, (\cantor ,\varphi )\in \mathbf{Q}\}\\
 & = & \mathcal{H}_{\mathbf{Q}}\times \mathcal{E}
\end{eqnarray*}
 In summary, we have the following commutative diagram:\[
\begin{array}{ccccc}
 \mathcal{H}_{\mathbf{Q}}\times \mathcal{E} & \xrightarrow{\beta } & \mathcal{B}_{\mathbf{Q}} & \subset  & \all _{\mathbf{Q}}\\
 \cap  &  & \cap  &  & \cap \\
 \mathcal{H}_{\mathbf{P}}\times \mathcal{E} & \xrightarrow{\beta } & \mathcal{B} & \subset  & \all _{\mathbf{P}}\end{array}
\]

We can now prove (1). Since $\mathcal{B}$ is residual in $\all _{\mathbf{P}}$,
we see that $\all _{\mathbf{Q}}$ is residual in $\all _{\mathbf{P}}$
if and only if $\mathcal{B}_{\mathbf{Q}}$ is residual in $\mathcal{B}$,
and since $\beta $ is a homeomorphism this is equivalent to $\beta ^{-1}(\mathcal{B}_{\mathbf{P}})=\mathcal{H}_{\mathbf{Q}}\times \mathcal{E}$
being residual in $\mathcal{H}_{\mathbf{P}}\times \mathcal{E}$. But
the latter happens if and only if $\mathcal{H}_{Q}$ is residual in
$\mathcal{H}_{\mathbf{P}}$. This completes the proof.

The proof of (2) is very similar. First note that $\mathcal{H}_{\mathbf{Q}}$
is dense in $\mathcal{H}_{\mathbf{P}}$ if and only if $\mathcal{H}_{\mathbf{Q}}\times \mathcal{E}$
is dense in $\mathcal{H}_{\mathbf{P}}\times \mathcal{E}$, which happens
if and only if $\beta (\mathcal{H}_{\mathbf{Q}}\times \mathcal{E})=\mathcal{B}_{\mathbf{Q}}$
is dense in $\mathcal{B}$. Since $\mathcal{B}$ is dense in $\all _{\mathbf{P}}$
this is equivalent to $\mathcal{B}_{\mathbf{Q}}$ being dense in $\all _{\mathbf{P}}$.

From this, one direction of (2) is immediate: if $\mathcal{H}_{\mathbf{Q}}$
is dense in $\mathcal{H}_{\mathbf{P}}$ then $\mathcal{B}_{\mathbf{Q}}$
is in $\all _{\mathbf{P}}$ and consequently $\all _{\mathbf{Q}}$
is too. Conversely, the remaining hypothesis about $\beta $ shows
that $\mathcal{B}_{\mathbf{Q}}$ is dense in $\all _{\mathbf{Q}}$,
since if $X\in \all _{\mathbf{Q}}$ then $X=\lim X_{n}$ with $X_{n}\in \mathcal{B}_{\mathbf{Q}}$;
so if $\all _{\mathbf{Q}}$ is dense in $\all _{\mathbf{P}}$, so
is $\mathcal{B}_{\mathbf{Q}}$.

This completes the proof, assuming the existence of $\mathcal{E}$
and $\beta $.
\end{proof}

\subsection{\label{sub:correspondence-construction}Construction of $\mathcal{E}$
and $\beta $.}

We turn to the details of the construction of $\mathcal{E}$ and $\beta $.
Write $\mathcal{D}$ for the space of all topological embeddings of
$\cantor $ in $Q$. This is a subspace of the space $\mathcal{C}=C(\cantor ,Q)$
of all continuous maps from $\cantor $ to $Q$, which carries the
usual metric \[
d_{\infty }(f,g)=\sup _{x\in \cantor }d(f(x),g(x))\]
for $f,g\in C(\cantor ,Q)$ (the symbol $d_{\infty }$ now represents
metrics on both $\mathcal{H}$ and $\mathcal{C}$; which is intended
will be clear from the context). This metric is complete on $\mathcal{C}$,
and one may verify that $\mathcal{D}$ is a $G_{\delta }$ subset
of $\mathcal{C}$ so it is a Polish space.

Our next lemma is a sharpening of a classical theorem of Kuratowski,
stating that the Cantor sets constitute a dense $G_{\delta }$ subset
of $2^{Q}$.

\begin{prop}
\label{prop:cantor-embedding} There exists a $G_{\delta }$ subset
$\mathcal{E}\subseteq \mathcal{D}$ such that the map $f\mapsto f(\cantor )$
from $\mathcal{E}$ to $\Cantor $ is an embedding, and its image
is dense.
\end{prop}
\begin{proof}
We construct a sequence $N_{i}\in \mathbb{N}$ and a family $(U_{\alpha })_{\alpha \in A}$
of open sets in $Q$ whose index set $A$ consists of finite words
\[
\alpha _{1}\ldots \alpha _{k}\in \{1,2,\ldots ,N_{1}\}\times \{1,2,\ldots ,N_{2}\}.\times \ldots \times \{1,2,\ldots ,N_{k}\}\]
 Write $\alpha _{1}\ldots \alpha _{n}\prec \alpha '_{1}\ldots \alpha '_{m}$
if $n\geq m$ and $\alpha _{1}\ldots \alpha _{m}=\alpha '_{1}\ldots \alpha '_{m}$.
Let $\emptyset $ denote the empty word, and write $A_{k}$ for all
words of length $k$ in $A$. We construct the $U_{\alpha }$ so that
they satisfy the following conditions: 
\begin{enumerate}
\item For each $k$ the family $\{U_{\alpha }\}_{\alpha \in A_{k}}$ is
pairwise disjoint, and if $\alpha \prec \alpha '$ then $U_{\alpha }\subseteq U_{\alpha '}$
(so the partially ordered set $(\{U_{\alpha }\},\subseteq )$ is isomorphic
to $(A,\prec )$) .
\item For each $k$, the union $\cup _{\alpha \in A_{k}}U_{\alpha }$ is
an open dense set in $Q$.
\item $\diam U_{\alpha }\leq 1/k$ for $\alpha \in A_{k}$ (This is the
condition that will determine the value of the $N_{i}$'s).
\end{enumerate}
Such a family can easily be constructed by recursion on $k$ by dividing
each $U_{\alpha }$, $\alpha \in A_{k}$ into finitely many small
disjoint open sets plus a meager remainder, obtaining the family $\{U_{\beta }\}_{\beta \in A_{k+1}}$.
We omit the details. Let \[
D=\bigcap _{k=1}^{\infty }\bigcup _{\alpha \in A_{k}}U_{\alpha }\]
which is dense $G_{\delta }$in $Q$. Define\[
\Cantor _{k}=\{C\subseteq \bigcup _{\alpha \in A_{k}}U_{\alpha }\, :\, C\textrm{ is a Cantor set}\}\]
 This set is open in $2^{Q}$. Finally, let \[
\Cantor =\bigcap _{k=1}^{\infty }\Cantor _{k}=\{C\subseteq D\, :\, C\textrm{ is a Cantor set}\}\]
 $\Cantor $ is a $G_{\delta }$ subset of $2^{Q}$ ; we next verify
that it is dense. Let $A\subseteq Q$ and $\varepsilon >0$, and let
$A_{\varepsilon }=\{x\in Q\, ;\, d(x,A)<\varepsilon \}$ and let $B=A_{\varepsilon }\cap D$.
Since $D$ is a dense $G_{\delta }$ in $Q$ and $A_{\varepsilon }$
is open we see that $B$ is a $G_{\delta }$ without isolated points
and is dense in $A_{\varepsilon }$. In particular $d(\overline{B},A)\leq \varepsilon $.
Since $B$ is relatively compact we can choose a finite set $x_{1},\ldots ,x_{N}\in B$
which is $\varepsilon $-dense in $B$. Since $B$ is a $G_{\delta }$
subset of a Polish space we may find Cantor subsets $C_{1},\ldots ,C_{N}\subseteq B$
with $C_{i}$ contained in the ball $B_{\varepsilon }(x_{i})$. Let
$C=\cup _{i=1}^{N}C_{i}$. Then $C\subseteq D$ is a Cantor set, so
$C\in \Cantor $, and $d(C,A)<3\varepsilon $. This establishes that
$\Cantor $is a dense $G_{\delta }$ in $2^{Q}$.

For each $C\in \Cantor $ we define a homeomorphism $\psi _{C}:C\rightarrow \cantor $
as follows. Assume without loss of generality that $\cantor \subseteq [0,1]$
is the standard middle-third realization of the Cantor set. Let $k$
be the first index such that for some $\alpha =\alpha _{1}\ldots \alpha _{k}\in A_{k}$
we have $C\subseteq U_{\alpha }$ and there are at least two distinct
indices $1\leq i,j\leq N_{k+1}$ such that $C\cap U_{\alpha i}\neq \emptyset $
and $C\cap U_{\alpha j}\neq \emptyset $. Let $r_{1},\ldots ,r_{n}$
be those indices such that $C\cap U_{\alpha r_{i}}\neq 0$; by our
assumption $n\geq 2$. Divide $\cantor $ into $n$ closed and open
sets $W_{1},\ldots ,W_{n}$ of diameter $<1/n$ in a manner depending
only on $k$ and $r_{1},\ldots ,r_{n}$ but not on $C$, and as a
first approximation prescribe $\psi _{C}$ maps $C\cap U_{\alpha r_{i}}$
to $W_{i}$. Now continue inductively to define a homeomorphism from
$C\cap U_{\alpha r_{i}}$ onto $W_{i}$. The map $\psi _{C}$ is defined
in the limit and is a homeomorphism from $C$ to $\cantor $.

We claim that $C\mapsto \psi _{C}^{-1}$ is an embedding of $\Cantor $
in $\mathcal{D}$. The inverse is continuous, since if $\psi _{C_{n}}^{-1}\rightarrow \psi _{C}^{-1}$
in $d_{\infty }$ then clearly $C_{n}\rightarrow C$ in the Hausdorff
metric. 

To see that $C\mapsto \psi _{C}^{-1}$ is continuous, let $C',C''\in \Cantor $
be close enough so that the first $k$ steps of the construction of
$\psi _{C'},\psi _{C''}$ agree. In particular there will be a pairwise
disjoint partition of $C'\cup C''$ by closed and open sets $V_{1},\ldots ,V_{n}$
each of diameter less than $1/k$, and a partition $W_{1},\ldots ,W_{n}$
of $\cantor $ into closed and open sets of diameter less than $1/2^{k}$,
such that both $\psi _{C'}$ maps $V_{i}\cap C'$ onto $W_{i}$ and
similarly $\psi _{C''}$ maps $V_{i}\cap C''$ onto $W_{i}$. Thus
$d(\psi _{C'}^{-1}(z),\psi _{C''}^{-1}(z))<1/k$ for every $z\in \cantor $,
so $d_{\infty }(\psi _{C'},\psi _{C''})<1/k$.

Thus the map $C\rightarrow \psi _{C}^{-1}$ is an embedding of $\Cantor _{\infty }$
in $\mathcal{D}$. Call its image $\mathcal{E}$. Since $\Cantor _{\infty }$
is a Polish space so is $\mathcal{E}$. This completes the proof.
\end{proof}
Given an embedding $f\in \mathcal{E}$ of $\cantor $ into $Q$ and
a homeomorphism $\varphi $ of $\cantor $, we can embed the system
$(\cantor ,\varphi )$ in $Q^{\mathbb{Z}}$ via the orbit picture
map $f_{\varphi }^{*}$ (section \ref{sub:orbit-iamge-map}). We define
$\beta :\mathcal{H}_{\mathbf{P}}\times \mathcal{E}\rightarrow \all $
in this way: $\beta (\varphi ,f)=f_{\varphi }^{*}(\cantor )$. Then
$(\beta (\varphi ,f),\sigma )$ is a system isomorphic to $(\cantor ,f)$.
Furthermore,

\begin{lem}
\label{lem:beta-is-an-embedding}The map $\beta :\mathcal{H}\times \mathcal{E}\rightarrow \all $,
which sends $(\varphi ,f)\in \mathcal{H}\times \mathcal{E}$ to the
orbit picture $f_{\varphi }^{*}(\cantor )$ of the systems $(\cantor ,\varphi )$,
is a topological embedding.
\end{lem}
\begin{proof}
Verification. 
\end{proof}
It remains to show that, if the Cantor systems are dense in $\all _{\textbf P}$,
then the image of $\mathcal{H}_{\textbf P}\times \mathcal{E}$ under
$\beta $ is dense in $\all _{\textbf P}$ . 

Recall that $\pi _{0}:Q^{\mathbb{Z}}\rightarrow Q$ is the projection
onto the $0$-th coordinate. $\pi _{0}$ is a continuous open and
closed map. A system $X\in \all $ is a \emph{graph} if there is a
homeomorphism $\varphi :\pi _{0}(X)\rightarrow \pi _{0}(X)$ such
that $X$ is the orbit picture of the system $(\pi _{0}(X),\varphi )$
under the inclusion map $i:\pi _{0}(X)\hookrightarrow Q$. Denote\[
\cgr =\{X\in \all \, :\, X\textrm{ is a graph}\}\]
We remark without proof that $\cgr $ is a dense $G_{\delta }$ subset
of $\all $ (but note that being a graph isn't an isomorphism invariant,
but rather a property of the embedding of the system in $\all $,
so it is not a dyamical property in our sense).

\begin{prop}
\label{lemma:cgr-denseness-in-zz}Let $\Cantor _{0}\subseteq \Cantor $
a dense set of Cantor sets in $Q$. Let $X\in \all $ be a Cantor
system. Then $X=\lim X_{n}$ for a sequence $X_{n}\in \cgr $ with
$\pi _{0}(X_{n})\in \Cantor _{0}$ and $(X_{n},\sigma )\cong (X,\sigma )$.
\end{prop}
\begin{proof}
Fix $\varepsilon >0$ and find $C\in \Cantor _{0}$ such that $d(C,\pi _{0}(X))<\varepsilon $
in $2^{Q}$; this can be done because $\Cantor _{0}$ is dense in
$2^{Q}$. Since $X$ is zero dimensional and $\pi _{0}:Q^{\mathbb{Z}}\rightarrow Q$
is open and closed, $\pi _{0}(X)$ is zero dimensional. Since $d(C,\pi _{0}(X))<\varepsilon $
and $X$ is a Cantor set, there exists a homeomorphism $\alpha :X\rightarrow C$
such that\[
d(x(0),\alpha x)<\varepsilon \]
for all $x\in X$. Define $T:C\rightarrow C$ by $T=\alpha \sigma \alpha ^{-1}$,
so $(C,T)\cong (X,\sigma )$. Let $i:C\rightarrow Q$ be the inclusion
map and define \[
Y=i_{T}^{*}(C)\]
the orbit picture of the system $(C,T)$ in $Q^{\mathbb{Z}}$. By
definition $Y\in \cgr $, and $(Y,\sigma )\cong (X,\sigma )$ via
the homeomorphism $j=i_{T}^{*}\circ \alpha :X\rightarrow Y$. We have
\[
(jx)(0)=(i^{*}\alpha x)(0)=\alpha x\]
so\[
d(x(0),(jx)(0))=d(x(0),\alpha x)<\varepsilon \]
 and more generally, since $\varphi $ commutes with $\sigma $ we
have\[
d(x(k),(jx)(k))=d((\sigma ^{k}x)(0),(j\sigma ^{k}x)(0))<\varepsilon \]
for all $x\in X$ and $k\in \mathbb{Z}$. Consequently, $d(X,Y)=d(X,jX)<\varepsilon $
in $\all $. 
\end{proof}
\begin{cor}
\label{cor:correspondence-dense-image}Under the hypotheses of the
correspondence theorem and with $\mathcal{E}$ as in lemma \ref{prop:cantor-embedding}
and $\beta $ as in lemma \ref{lem:beta-is-an-embedding}, the image
of $\mathcal{H}_{\textbf P}\times \mathcal{E}$ under $\beta $ is
dense in $\all _{\textbf P}$
\end{cor}
\begin{proof}
Write $\mathcal{B}$ for the image of $\mathcal{H}_{\textbf P}\times \mathcal{E}$
under $\beta $ and let $\Cantor _{0}=\{f(\cantor )\, :\, f\in \mathcal{E}\}\subseteq \Cantor $.
By assumption there is a dense set $\Cantor '$ of Cantor systems
in $\all _{\textbf P}$, so it suffices to show that the $\mathcal{B}$
is dense in $\Cantor '$. If $X\in \Cantor '$ is a cantor system,
then by the lemma above $X=\lim X_{n}$ for $X_{n}\in \cgr $, $\pi _{0}(X_{n})\in \Cantor _{0}$
and $(X,\sigma )\cong (X_{n},\sigma )$. The last fact implies that
$(X_{n},\sigma )\in \all _{\textbf P}$. Let $\psi _{n}:\pi _{0}(X_{n})\rightarrow \pi _{0}(X_{n})$
be a homeomorphism such that $X_{n}$ is the orbit picture of $\psi _{n}$,
and let $f\in \mathcal{E}$ such that $f:\cantor \rightarrow \pi _{0}(X_{n})$
is a homeomorphism. Define $\varphi _{n}\in \mathcal{H}$ by $\varphi _{n}=f^{-1}\psi _{n}f$.
Since $(\cantor ,\varphi _{n})\cong (X_{n},\sigma )$ we have $\varphi _{n}\in \mathcal{H}_{\textbf P}$,
and $\beta (\varphi _{n},f)=X_{n}\in \mathcal{B}$.
\end{proof}
This completes the proof of the correspondence theorem. 

Let us show that the hypotheses of the correspondence theorem are
satisfied for the class $\textbf P$ of all systems, i.e. that the
genericity status of a dynamical property $\textbf Q$ is the same
in $\mathcal{H}$ and $\all $. Conditions (a) and (b) are satisfied
trivially; (c) follows from:

\begin{prop}
\label{pro:Cantor-systems-are-dense}The Cantor systems are dense
in $\all $.
\end{prop}
\begin{proof}
Let $X\in \all $ and $\varepsilon >0$, and write $D=\pi _{0}(X)$.
Let $C\in 2^{Q}$ be a Cantor set with $d(C,D)<\varepsilon $. Let
\[
X'=\{x'\in C^{\mathbb{Z}}\, :\, \exists x\in X\textrm{ s.t. }d(x(i),x'(i))\leq \varepsilon \textrm{ for all }i\in \mathbb{Z}\}\]
This set is seen to be nonempty, closed and shift-invariant, so $X'\in \all $.
Clearly $d(X,X')<\varepsilon $. We also have $X'\subseteq C^{\mathbb{Z}}$
so $X'$ is zero dimensional, and it is easily seen to have no isolated
points, so $X'$ is a Cantor system.
\end{proof}
\begin{cor}
\label{cor:Kekhris-Rosendal-generic-in-S}The Kechris-Rosendal system
is generic in $\all $.
\end{cor}

\subsection{\label{sub:Zero-One-Laws}Zero-One Laws}

We turn now to the zero-one laws. Since there are sets which are neither
meagre nor residual, one would expect that there are dynamical properties
which are neither generic nor exotic in $\all $. In \cite{GK-01law}
Glasner and King proved a zero-one law asserting that in the automorphism
group $\mathcal{A}$ of a Lebesgue space every Baire measurable dynamical
property is either generic or exotic. We show next that a similar
situation holds for certain subspaces of $\mathcal{H}$. 

A subset $A$ of a complete metric space $X$ is \emph{Baire measurable}
if it belongs to the completion of the Borel $\sigma $-algebra of
$X$ with respect to the $\sigma $-ideal of meagre sets. Equivalently,
it is Baire measurable if it can be written as $U\Delta M$ where
$U$ is open and $M$ meagre. We follow the convention from \cite{GK-01law}
and use the above terminology instead of the usual {}``sets with
the property of Baire''. All Borel sets, and more generally all analytic
sets, are Baire measurable. For more information see Oxtoby \cite{Ox}.

\begin{thm}
\label{thm:general-zero-one-law}Let $\mathbf{P}$ be a dynamical
property such that $\mathcal{H}_{\mathbf{P}}$ is a $G_{\delta }$
and there is a Cantor system $(C,\varphi )$ whose the isomorphism
class is dense in $\mathcal{H}_{\mathbf{P}}$. Then for every dynamical
property $\mathbf{Q}\subseteq \mathbf{P}$ with $\mathcal{H}_{\mathbf{Q}}$
a Baire set, either $\mathbf{Q}$ is generic in $\mathcal{H}_{\mathbf{P}}$
or it is exotic there.
\end{thm}
This is an immediate corollary of the following:

\begin{thm}
\label{thm:invariant-sets-of-group-actions}(Glasner and King, \cite{GK-01law})
Suppose a group $\Gamma $ acts by homeomorphisms on a Polish space
$X$. If the action is transitive (i.e. there is a $x\in X$ such
that $\overline{\Gamma x}=X$) then every Baire measurable subset
of $X$ which is invariant under the action of $\Gamma $ is either
of meagre or residual.
\end{thm}
The proof is not complicated and can be found in \cite{GK-01law}.

\begin{proof}
(of theorem \ref{thm:general-zero-one-law}) Let $\mathcal{H}$ act
on itself by conjugation, and note that both $\mathcal{H}_{\mathbf{P}}$
and $\mathcal{H}_{\mathbf{Q}}$ are invariant under this action. Furthermore,
by assumption $\mathcal{H}_{\mathbf{P}}$ is Polish, the $\mathcal{H}$-
orbit of $\varphi $ is dense in $\mathcal{H}_{\mathbf{P}}$, and
$\mathcal{H}_{\mathbf{Q}}$ is a $\mathcal{H}$-invariant Baire set;
so by theorem \ref{thm:invariant-sets-of-group-actions} it is either
residual or meagre. 
\end{proof}
Theorem \ref{thm:zero-one-law} now follows from part \ref{thm:totally-transitive-systems-rokhlin}
of theorem \ref{thm:totally-transitive-systems} and the theorem above.

\section{\label{sec:symbolic-approximation}Symbolic approximation}

\subsection{\label{sub:Symbolic-Systems}Symbolic Systems }

Recall that $\pi _{0}:Q^{\mathbb{Z}}\rightarrow Q$ is the projection
onto the $0$-th coordinate. A system $X\in \all $ is \emph{symbolic}
if its cross-section $\pi _{0}(X)$ is finite. In this case $\pi _{0}(X)$
is also called the \emph{alphabet} of $X$. 

Write $Q^{*}$ for the set of all finite words over $Q$. We denote
words by the letters $a,b,c$ or $u,v,w$. The $i$-th letter of a
word $a$ is $a(i)$. If $a=a(1)a(2)\ldots a(k)$ then $k$ is the
length of $a$ and is denoted by $\ell (a)$. The concatenation of
words $a,b\in Q^{*}$ is written $ab$. For $a,b\in Q^{*}$ we say
that $a$ is a subword of $b$ at index $i\leq \ell (b)-\ell (a)+1$
if $a(j)=b(i+j)$ for $j=1,\ldots ,\ell (a)$. The index $i$ is called
the \emph{alignment} of $a$ in $b$. If such an $i$ exists, we say
that $a$ appears in $b$, or is a subword of $b$.

We measure the distance between two words $a,b$ with the same length
(finite or infinite) by \[
d_{\infty }(a,b)=\sup _{i}d(a(i),b(i))\]
 (the symbol $d_{\infty }$ was defined already in section \ref{sec:correspondence};
the new meaning can be distinguished from the old one by its context).
For $\varepsilon >0$ we say the the word $a$ is an $\varepsilon $-subword
in $b$ if there is a word $a'$ appearing in $b$ with $d_{\infty }(a,a')<\varepsilon $.
Such a subword $a'$ of $b$ is called an $\varepsilon $-appearance
of $a$ in $b$. 

Note that $x,y\in Q^{\mathbb{Z}}$ are $\varepsilon $-close, i.e.
$d(x,y)<\varepsilon $ with respect to the metric on $Q^{\mathbb{Z}}$
introduced in section \ref{sec:Definitions}, if and only if for $n=\left\lceil 1/\varepsilon \right\rceil $
the words $a=x(-n)\ldots x(n)$ and $b=y(-n)\ldots y(n)$ satisfy
$d_{\infty }(a,b)<\varepsilon $.

Every $X\in \all $ is a set of bi-infinite sequences, which we think
of as bi-infinite words over $Q$. We say that a finite word $a=a(1)\ldots a(n)$
appears (or $\varepsilon $-appears) in $X$ if there is an $x\in X$
such that $a$ appears (or $\varepsilon $-appears) in $x$. By shift
invariance of $X$, if $a$ is a subword ($\varepsilon $-subword)
of $X$ then there are appearances ($\varepsilon $-appearances) of
$a$ in $X$ with every alignment.

Using this terminology, for $X,Y\in \all $ we have $d(X,Y)\leq \varepsilon $
if and only if whenever $a$ is a subword of $X$ and $\ell (a)\leq 1+2/\varepsilon $
then $a$ is an $\varepsilon $-subword of $Y$, and similarly with
the roles of $X,Y$ reversed.

\subsection{\label{sub:Approximation-by-symbolic}Approximation by symbolic systems}

A finite set of finite words $L\subseteq Q^{*}$ is called a \emph{Language.}
A system $X\in \all $ is said to be \emph{constructed from} $L$
if every word in $X$ is a bi-infinite concatenation of words from
$L$ and every word from $L$ appears in $X$. 

Given a system $X\in \all $ our first goal is to find workable conditions
under which a system $Y$ constructed from a language $L$ is close
to $X$ in $\all $. The basic idea will be to work with languages
$L$ which are made up of words which appear (or $\varepsilon $-appear)
in $X$. Thus a concatenation of words from $L$ will look locally
like a subword (or $\varepsilon $-subword) of $X$, provided we take
care not to splice together subwords of $X$ which don't {}``fit''.
We will also want $L$ to be large enough that it contain words representing
all subwords of $X$, up to some degree of accuracy. 

The following notion is central to making the above precise. Let $(X,T)$
be a dynamical system and $\varepsilon >0$. A finite, infinite or
bi-infinite sequence $\overline{x}=(x_{i})$ is called an \emph{$\varepsilon $-pseudo-orbit}
if $d(x_{i+1},Tx_{i})<\varepsilon $ for all $i$. A finite $\varepsilon $-pseudo
orbit $x_{1}\ldots x_{N}$ is an $\varepsilon $-\emph{pseudo-period}
if in addition $d(Tx_{N},x_{1})<\varepsilon $. 

For $X\in \all $ and a pseudo-orbit $\overline{x}=(x_{i})$ in $X$,
note that each $x_{i}$ is itself a bi-infinite sequence over $Q$,
$x_{i}=(\ldots ,x_{i}(-1),x_{i}(0),x_{i}(1)\ldots )$. We adopt the
convention that sequences of points in $Q^{\mathbb{Z}}$ (finite,
one-sided infinite or bi-infinite) are always written using the bar
notation $\overline{x}$. The member points of such a sequence are
written as $x_{i}$, with added superscripts if necessary; this $x_{i}(j)\in Q$
is the $j$-th coordinate of the $i$-th point $x_{i}\in Q^{\mathbb{Z}}$
of the sequence $\overline{x}$.

Let $X\in \all $. Given a bi-infinite sequence $\overline{x}=(\ldots ,x_{-1},x_{0},x_{1},\ldots )$
of points in $X$ we define $\theta (\overline{x})\in Q^{\mathbb{Z}}$
by $\theta (\overline{x})(i)=x_{i}(0)$, so \[
\theta (\ldots ,x_{-1},x_{0},x_{1}\ldots )=(\ldots ,x_{-1}(0),x_{0}(0),x_{1}(0),\ldots )\]
 We define $\theta $ similarly on finite sequences of points from
$X$, so $\theta (x_{1},\ldots ,x_{N})$ is a finite word $a$ of
length $N$ with $a(i)=x_{i}(0)$.

\begin{lem}
\label{lem:approximation-bound}If $\overline{x}=(\ldots x_{-1},x_{0},x_{1}\ldots )$
is a bi-infinite $\varepsilon $-pseudo-orbit in $(Q^{\mathbb{Z}},\sigma )$
and $\varepsilon <1$, then $d(\theta (\overline{x}),x_{0})\leq \sqrt{\varepsilon }$
\end{lem}
\begin{proof}
Suppose $k\leq 1/\sqrt{\varepsilon }$. It suffices to show that $d(x_{0}(i),\theta (\overline{x})(i))<\sqrt{\varepsilon }$
for $|i|\leq k$. But\begin{eqnarray*}
d(x_{0}(i),\theta (\overline{x})(i)) & = & d(x_{0}(i),x_{i}(0))\\
 & \leq  & \sum _{j=0}^{i-1}d(x_{j}(i-j),x_{j+1}(i-(j+1)))
\end{eqnarray*}
 since $d(\sigma x_{j},x_{j+1})<\varepsilon $, for every $m\leq 1/\varepsilon $
it holds that $d(\sigma x_{j}(m),x_{j+1}(m))=d(x_{j}(m+1),x_{j+1}(m))<\varepsilon $
and since $k\leq 1/\sqrt{\varepsilon }$ this is certainly true for
$m=i-j$ as $j$ ranges from $0$ to $i$. Therefore each summand
is less than $\varepsilon $ so\[
d(x_{0}(i),\theta (\overline{x})(i))<i\varepsilon \leq \sqrt{\varepsilon }\]
 as claimed. 
\end{proof}
Let $X\in \all $ and for $i=1,\ldots ,M$ suppose that $\overline{x}^{i}=x_{1}^{i},\ldots ,x_{N(i)}^{i}$
are $\varepsilon $-pseudo-orbits with $\varepsilon <1$. Let $a^{i}=\theta (\overline{x}^{i})$
and $L=\{a^{1},\ldots ,a^{M}\}$. Suppose $Y\neq \emptyset $ is a
system constructed from $L$ and that $a^{i}a^{j}$ appears in $Y$
only if $d(\sigma (x_{N(i)}^{i}),x_{1}^{j})<\varepsilon $. Every
$y\in Y$ can be written as a concatenation of the form $y=\ldots a^{k(-1)}a^{k(0)}a^{k(1)}\ldots $,
and then $y=\theta (\overline{y})$ where $\overline{y}$ is a bi-infinite
sequence of points from $X$ given by the concatenation $\overline{y}=\ldots \overline{x}^{k(-1)}\overline{x}^{k(0)}\overline{x}^{k(1)}\ldots $
aligned in the obvious way. By the assumption that $a^{i}a^{j}$ appears
in $Y$ only if $d(\sigma (x_{N(i)}^{i}),x_{1}^{j})<\varepsilon $
we $\overline{y}$ is an $\varepsilon $-pseudo-orbit. The lemma now
implies that every $y\in Y$ is $\sqrt{\varepsilon }$-close to one
of the $x_{j}^{i}$, so every $y\in Y$ is $\sqrt{\varepsilon }$-close
to a point in $X$. 

Now suppose in addition that the union of the collection $\{x_{j}^{i}\, :\, 1\leq i\leq M\, ,\, 1\leq j\leq N(i)\}\subseteq X$
of points making up the pseudo-orbits $\overline{x}^{1},\ldots ,\overline{x}^{M}$
is $\varepsilon $-dense in $X$, that is, for every $x\in X$ we
have $d(x,x_{j}^{i})<\varepsilon $ for some $i,j$. Recall that we
are assuming that $Y$ was constructed from $L$, so every $a^{i}$
appears in $Y$. We claim that this ensures that $d(x,Y)<2\sqrt{\varepsilon }$
for every $x\in X$. Indeed, let $x\in X$ be arbitrary. By assumption
there are $i,k$ such that $d(y,x_{i}^{k})<\varepsilon $. Therefore
there is a $y\in Y$ and an $\varepsilon $-pseudo-orbit $\overline{y}=(\ldots ,x_{0},x_{1},\ldots )$
with $x_{0}=x_{i}^{k}$ and $y=\theta (\overline{y})$ follows from
the denseness of the $x_{j}^{i}$; . By the previous lemma, $d(\theta (\overline{y}),x_{i}^{k})=d(\theta (\overline{y}),x_{0})\leq \sqrt{\varepsilon }$
so \[
d(x,y)\leq d(x,x_{i}^{k})+d(x_{i}^{k},y)=d(x,x_{i}^{k})+d(\theta (\overline{y}),x_{i}^{k})\leq \varepsilon +\sqrt{\varepsilon }\leq 2\sqrt{\varepsilon }\]
 We have proved the following:

\begin{prop}
\label{prop:symbolic-approximation}Let $X\in \all $ and suppose
that for $i=1,\ldots ,M$ we are given $\varepsilon $-pseudo-orbits
$\overline{x}^{i}=x_{1}^{i},\ldots ,x_{N(i)}^{i}$ in $X$ such that
their union is $\varepsilon $-dense in $X$. Write $a^{i}=\theta (\overline{x}^{i})$
and $L=\{a^{1},\ldots ,a^{M}\}$ and let $Y$ be any system constructed
from $L$ under the restriction that $a^{i}a^{j}$ appears in $Y$
only if $d(\sigma (x_{N(i)}^{i}),x_{1}^{j})<\varepsilon $. Then $d(X,Y)<2\sqrt{\varepsilon }$. 
\end{prop}
A special and very useful case of this is:

\begin{cor}
\label{cor:symbolic-approximation}Let $X\in \all $ and suppose that
for $i=1,\ldots ,M$ we are given $\varepsilon $-pseudo-periods $\overline{x}^{i}=x_{1}^{i},\ldots ,x_{N(i)}^{i}$
in $X$ with the same initial point, i.e. $x_{1}^{1}=x_{1}^{2}=\ldots =x_{1}^{M}$.
Assume that the union of the pseudo-periods is $\varepsilon $-dense
in $X$. Write $a^{i}=\theta (\overline{x}^{i})$ and $L=\{a^{1},\ldots ,a^{M}\}$
and let $Y$ be any system constructed from $L$. Then $d(X,Y)<2\sqrt{\varepsilon }$. 
\end{cor}

\subsection{\label{sub:Projection-into-symbolic-systems}Projection into symbolic
systems}

So far we have seen how to construct systems $Y$ close to $X$, provided
we have a good language to work with. We now would like this $Y$
to have the property that any system $Z$ sufficiently close to $Y$
inherits some of $Y$'s structure. This is our next task.

Let $Y\in \all $ and $a\in Q$. If there is in $\pi _{0}(Y)$ a unique
point closest to $a$, denote it by $\tau _{Y}(a)$. This defines
a partial map $\tau _{Y}:Q\rightarrow \pi _{0}(Y)\subseteq Q$, which
is defined on an open subset of $Q$ and is continuous there. Extend
this to a partial map $\tau _{Y}:Q^{\mathbb{Z}}\rightarrow Q^{\mathbb{Z}}$
by $\tau _{Y}(z)(i)=\tau _{Y}(z(i))$. When $\tau _{Y}$ is defined
on a point $z\in Q^{\mathbb{Z}}$ it is defined also on $\sigma ^{k}z$
for all $k$ and commutes with $\sigma $ in the sense that $\sigma (\tau _{Y}z)=\tau _{Y}(\sigma z)$;
also, $\tau _{Y}$ is continuous where it is defined (though its domain
may not be open). In particular, if $Z\in \all $ is such that $\tau _{Y}$
is defined on every $z\in Z$, we see that $\tau _{Y}$ is a factor
map from $Z$ onto some subsystem of $Y$, and that if $Y$ is symbolic
then $\tau _{Y}(Z)\subseteq \pi _{0}(Y)^{\mathbb{Z}}$, and so $\tau _{Y}(Z)$
is symbolic over the same alphabet as $Y$.

\begin{lem}
\label{lem:symbolic-retract}If $Y\in \all $ is symbolic then $\tau _{Y}$
is defined on every $Z\in \all $ sufficiently close to $Y$. 
\end{lem}
\begin{proof}
Let $Y,Z\in \all $ with $Y$ symbolic. Define \begin{eqnarray*}
\rho _{Y} & = & \min \{d(a,b)\, :\, a,b\in \pi _{0}(Y)\textrm{ and }a\neq b\}
\end{eqnarray*}
 $\pi _{0}(Y)$ is finite, hence the minimum in the definition of
$\rho _{Y}$ exists. Suppose that $d(Z,Y)<\frac{1}{2}\rho _{Y}$.
If $z\in Z$ then there is a $y\in Y$ with $d(z,y)<\frac{1}{2}\rho _{Y}$,
so $d(z(0),y(0))<\frac{1}{2}\rho _{Y}$. Thus clearly $\tau _{Y}$
is defined on $z(0)$. Since $z\in Z$ was arbitrary and using the
shift invariance of $Z$ we conclude that $\tau _{Y}$ is is defined
on every $z\in Z$. 
\end{proof}
\label{def:unique-parsing}We say that a finite language $L\subseteq Q^{*}$
has the \emph{unique parsing property} \emph{with window size} $N$
if every word $a\in Q^{*}$ of length at least $N$ has at most one
parsing $a=uv_{1}\ldots v_{m}w$ such that (a) $v_{1},\ldots ,v_{m}\in L$,
and (b) $u'u,ww'\in L$ for some $u',w'\in Q^{*}$. 

If $L$ is a language such that distinct words in $L$ contain distinct
letters, then $L$ has the unique parsing property. Since $Q$ has
no isolated points, we can always perturb the letters of words of
a finite alphabet $L$ by an arbitrarily small amount to make them
distinct, and achieve unique parsing in this way. 

Another way to get unique parsing is to add prefixes. Suppose, for
instance, that $u,v$ are two words. Let $w$ be a third word and
$n$ such that $w^{n}=w\ldots w$ ($n$ times) does not appear as
a subword of $u,v$. Then setting $u'=w^{n}u$ and $v'=w^{n}u$, we
see that $L=\{u',v'\}$ has the unique parsing property with window
size $N=n+\max \{\ell (u'),\ell (v')\}$. 

There are many other ways to get unique parsing. In the sequel we
will not spell out the details of this.

Suppose $Y\in \all $ is a symbolic system constructed from a language
$L$ which has the unique parsing property with window size $N$.
Let $Z$ be close enough to $Y$ that for every $z\in Z$ there is
a $y\in Y$ such that $d(z(i),y(i))<\rho _{Y}/2$ for every $|i|\leq N$,
with $\rho _{Y}$ as in the proof of the previous lemma. It follows
that for every $z\in Z$ the central $(2N+1)$-long subword of $\tau _{Y}(z)$
equals $y(-N),\ldots ,y(N)$ for some $y\in Y$. Given $z\in Z$ let
$a=\tau _{Y}(z)$ and $a_{i}=a(i-N)\ldots a(i+N)$. Each $a_{i}$
has a unique parsing as in definition \ref{def:unique-parsing}, and
by uniqueness the parsing of the word $a(i-N+1),\ldots ,a(i+N)$ induced
by the parsings of $a_{i}$ and of $a_{i+1}$ must agree. This means
that we can merge all the parsings of the $a_{i}$'s and obtain a
parsing of $a$ into words from $L$. This proves:

\begin{lem}
\label{lem:retract-embedding}Let $Y$ be constructed from a language
$L$ with the unique parsing property. If $Z$ is close enough to
$Y$ (in a manner depending only on $L$), then $\tau _{Y}$ is defined
on $Z$ and $\tau _{Y}(Z)$ is also constructed from $L$. 
\end{lem}
There is one last important property of $\tau _{Y}$ that we will
use. For $x,y\in Q^{\mathbb{Z}}$ let \[
d_{\infty }(x,y)=\sup _{k\in \mathbb{Z}}d(x(k),y(k))\]
 denote the uniform distance between $x$ and $y$. 

\begin{lem}
\label{lem:uniform-perturbation}Let $Y\in \all $ be a symbolic system
and $\varepsilon >0$. Then for every system $Z$ sufficiently close
to $Y$, the projection $\tau _{Y}:Z\rightarrow Y$ is defined and
displaces points by at most $\varepsilon $, that is, $d_{\infty }(z,\tau _{Y}(z))<\varepsilon $
for all $z\in Z$.
\end{lem}
\begin{proof}
It is easy to check that if $d(Y,Z)<\varepsilon $ then $d(y(0),\pi _{0}(Z))<\varepsilon $
for every $y\in Y$, implying $d(z(0),\tau _{Y}(z)(0))<\varepsilon $.
Of course there is nothing special about the index $0$, so assuming
that $\tau _{Y}$ is defined on $Z$ (which will be true for $Z$
close enough to $Y$) we have $d_{\infty }(z,\tau _{Y}(z))<\varepsilon $
for all $z\in Z$.
\end{proof}

\section{\label{sec:transitivity-and-periodic-approximation}The Space of
Transitive Systems}

In this section we prove theorem \ref{thm:transitive-systems}. We
begin with some definitions.

\subsection{\label{sub:Transitivity}Transitivity}

A dynamical system $(X,T)$ is \emph{bi-transitive} if there is a
point $x\in X$ whose full orbit $\{T^{n}x\}_{n\in \mathbb{Z}}$ is
dense in $X$; such a point is called a bi-transitive point. The system
is \emph{forward transitive} if there is a point $x\in X$ whose forward
orbit $\{T^{n}x\}_{n\in \mathbb{N}}$ is dense; such a point is called
a \emph{forward transitive point}. A point $x\in X$ is \emph{forward
recurrent} if there is a sequence $n(k)\in \mathbb{N}$ with $n(k)\rightarrow \infty $
and $T^{n(k)}x\rightarrow x$. 

For brevity, we say that a system is \emph{transitive} if it is forward
transitive, and denote this class of systems by $\trans $. 

The notions of bi-transitivity and forward-transitivity are distinct
in general, but for systems without isolated points (and in particular
the Cantor set), they coincide. The following sequence of lemmas establishes
this, along with some other well known facts we will use. For completeness
we provide proofs; or see e.g. \cite{Wal81}.

\begin{lem}
\label{lem:implications-of-forward-transitivity}Let $(X,T)$ be forward
transitive. Then every forward transitive point is forward recurrent,
and either $X$ is perfect or else $X$ is finite and $T$ is a cyclical
permutation of $X$.
\end{lem}
\begin{proof}
Let $x\in X$ be forward transitive. If $T^{n}x=x$ for some $n>0$
then the forward orbit of $x$, and hence $X$, is $\{x,Tx,\ldots ,T^{n-1}x\}$,
and $T$ permutes this set cyclically. We also see that $T^{nk}x\rightarrow x$
as $k\rightarrow \infty $, so $x$ is forward recurrent. Otherwise
$T^{n}x\neq x$ for all $n\geq 1$, so $T^{n}x\neq T^{-1}x$ for all
$n\geq 0$. Since $T^{-1}x$ is in the closure of $\{T^{n}x\}_{n\geq 0}$
there is a sequence $n_{k}\rightarrow \infty $ with $T^{n_{k}}x\rightarrow T^{-1}x$,
so $T^{n_{k}+1}x\rightarrow x$ , showing that $x$ is forward recurrent.
Since $T^{n_{k}+1}x\neq x$, this shows that $x$ is not isolated,
and, since $T$ is a homeomorphism, neither are any of the points
$T^{n}x$; since $\{T^{n}x\}_{n\in \mathbb{N}}$ is dense, $X$ has
no isolated points.
\end{proof}
\begin{lem}
\label{lem:charcterization-of-forward-transitivity}A dynamical system
$(X,T)$ is transitive if and only if for every nonempty nonempty
open sets $U,V\subseteq X$ there is an $n>0$ with $U\cap T^{-n}V\neq \emptyset $.
\end{lem}
\begin{proof}
Suppose $x$ is forward transitive. We may assume $X$ is perfect,
since otherwise $X$ consists of a single orbit and the conclusion
holds trivially. Hence $U,V$ are infinite, so the forward orbit $\{T^{n}x\}_{n\in \mathbb{N}}$
visits each of them infinitely often (otherwise $U\setminus \{T^{n}x\}_{n=1}^{\infty }$
would be a nonempty open set not visited by $x$'s forward orbit;
and similarly for $V$). We can therefore find integers $0<i<j$ with
$T^{i}x\in V$ and $T^{j}x\in U$. Hence $T^{j}x\in U\cap T^{i-j}V$,
so $U\cap T^{i-j}V\neq \emptyset $.

Conversely, let $\{U_{i}\}_{i\in I}$ be a countable basis for $X$,
and $U_{i}\neq \emptyset $. It is clear that the set of forward transitive
points of $X$ is precisely $\cap _{i\in I}\cup _{n\in \mathbb{N}}T^{-n}U_{i}$.
By assumption $\cup _{n\in \mathbb{N}}T^{-n}U_{i}$ is dense in $X$
for all $i\in I$. Hence $\cap _{i\in I}\cup _{n\in \mathbb{N}}T^{-n}U_{i}$
is a dense $G_{\delta }$, and in particular nonempty.
\end{proof}
\begin{lem}
\label{lem:transitivity-in-perfect-spaces}If $X$ is has no isolated
points, then $(X,T)$ is bi-transitive if and only if it is forward
transitive.
\end{lem}
\begin{proof}
One direction is trivial. In the other, suppose then that $(X,T)$
is transitive and $X$ is perfect. Let $U,V\neq \emptyset $ be open
sets; we must show that $U\cap T^{-n}V\neq \emptyset $ for some $n>0$.
Since $X$ is perfect $U,V$ are infinite. The orbit of $x$ must
visit $U$ infinitely many times, since otherwise $U\setminus \{T^{n}x\, ;\, n\in \mathbb{Z}\}$
is open, nonempty and not visited by the orbit of $x$. The same holds
for $V$. Therefore there are integers $i<j$ such that $T^{i}x\in V$
and $T^{j}x\in U$; so $T^{j}x\in U\cap T^{-n}V\neq \emptyset $ for
$n=j-i>0$.
\end{proof}
Since the Cantor set has no isolated points, the spaces of bi-transitive
and of forward transitive systems in $\mathcal{H}=\homeo (\cantor )$
coincide. This equality does not hold in $\all $. To see this, let
$0,1\in Q$ denote distinct elements and let $x=(\ldots ,0,0,0,1,1,1,\ldots )\in Q^{\mathbb{Z}}$
and $X$ the closure of the full orbit of $X$. Then $X$ consists
of the fixed points $(\ldots ,0,0,0,\ldots )$ and $(\ldots ,1,1,1,\ldots )$
and a single orbit spiralling between them, all of whose points are
bi-transitive but not forward transitive. 

We note without proof that the space of transitive systems in $\all $,
which is a Polish subset, generically contains just one system, and
this systems is not forward transitive. We henceforth concentrate
on the space of forward transitive systems.

\subsection{\label{sub:the-spaces-of-transitive-systems}The Spaces of Transitive
Systems are $G_{\delta }$'s.}

In order to apply the correspondence theorem to the spaces of forward-transitive
systems we must first establish that these spaces are $G_{\delta }$'s.
Recall that the property of being (forward) transitive is denoted
by $\trans $. 

\begin{prop}
\label{pro:trans-is-Gdelta-in-G}$\mathcal{H}_{\trans }$ is a $G_{\delta }$
subset of $\mathcal{H}$.
\end{prop}
\begin{proof}
Let $\{U_{i}\}_{i\in I}$ be a countable basis for the Cantor set,
with $U_{i}\neq \emptyset $. For each $i,j\in I$ and $n\in \mathbb{N}$,
the set of $\varphi \in \homeo (\cantor )$ such that\[
\varphi ^{-n}(U_{i})\cap U_{j}\neq \emptyset \]
is an open set in $\mathcal{H}$, and thus for each $i,j$ the set
of $\varphi \in \mathcal{H}$ so that $\varphi ^{-n}(U_{i})\cap U_{j}\neq \emptyset $
for some $n\in \mathbb{N}$ is open. Intersecting over $i,j\in I$
we get a $G_{\delta }$, and this is precisely the set of transitive
homeomorphisms by \ref{lem:charcterization-of-forward-transitivity}.
\end{proof}
The proof for $\all $ is similar:

\begin{prop}
\label{thm:trans-is-Gdelta-in-S}$\all _{\trans }$ is a $G_{\delta }$
subset of $\all $
\end{prop}
\begin{proof}
Let $\{U_{i}\}_{i\in I}$ be a countable basis for the topology of
$Q^{\mathbb{Z}}$. A system $X\in \all $ is transitive if and only
if for every $i,j\in I$ \begin{equation}
X\cap U_{i}\neq \emptyset \textrm{ and }X\cap U_{j}\neq \emptyset \; \Rightarrow \; \exists k>0\, :\, X\cap U_{i}\cap \sigma ^{-k}U_{j}\neq \emptyset \label{eq:transitivity-condition}\end{equation}
Thus it suffices to show that for fixed $i,j\in I$, the systems that
have this property with respect to $U_{i},U_{j}$ are a $G_{\delta }$
set.

Fix $i,j\in I$. The set of systems which fail to intersect one of
$U_{i}$ or $U_{j}$ is closed in the Hausdorff metric, and hence
is a $G_{\delta }$ in $\all $.

Let $\mathcal{U}\subseteq \all $ denote the class of system which
intersect both $U_{i}$ and $U_{j}$; by the above this set is open
in $\all $. For each $k>0$, the set $\mathcal{U}_{k}$ of systems
$X\in \mathcal{U}$ such that $X\cap U_{i}\cap \sigma ^{-k}U_{j}\neq \emptyset $
is open. The set of systems which satisfy condition (\ref{eq:transitivity-condition})
is $\cup _{k>0}\mathcal{U}_{k}$, which is open (and of course a $G_{\delta }$).

We have shown that for fixed $i,j\in I$, the set of systems in $\all $
which satisfy the condition (\ref{eq:transitivity-condition}) above
consists of the union of two $G_{\delta }$ sets, and is hence itself
a $G_{\delta }$ (in general, the union of finitely many $G_{\delta }$'s
is a $G_{\delta }$).
\end{proof}
We remark that the spaces $\mathcal{H}_{\trans }$ and $\all _{\trans }$
are not closed. To see that $\all _{\trans }$ is not closed, let
$0,1\in Q$ be distinct points, let $0_{n},1_{n}$ be sequences of
$n$ repetitions of $0,1$ respectively, and let $x_{n}$ be the periodic
sequence $x_{n}=\ldots 0_{n}\, 1_{n}\, 0_{n}\, 1_{n}\ldots $. The
orbit closure $X_{n}$ of $x_{n}$ is transitive (it consists of a
single periodic orbit) but $X_{n}\rightarrow X$ where $X$ is the
simple system generated by the points $(\ldots 0,0,1,1\ldots )$ and
$(\ldots 1,1,0,0\ldots )$ and this systems is not transitive. One
can easily construct examples of this sort where the $X_{n}$ are
infinite. A similar construction may be carried out in $\mathcal{H}$
(see \cite{BDK}).

\subsection{\label{sub:Finite-Cycles}Finite Cycles}

The simplest transitive systems are those which are a cyclical permutation
of a finite set. We call such systems \emph{finite cycles}. A finite
cycle with $n$ points is isomorphic to the finite group $\mathbb{Z}/n\mathbb{Z}$
with the map $i\mapsto i+1\bmod n$. We denote this system by $C_{n}$
and use the same symbol to denote the group $\mathbb{Z}/n\mathbb{Z}$. 

Note that a system $Y\in \all $ is a finite cycle if and only if
it consists of the translates of a single periodic $Q$-sequence.

It is well known that in the space $\mathcal{M}$ of invariant probability
measures on $Q^{\mathbb{Z}}$, the measures supported on periodic
orbits are dense. Below we prove a topological analogue of this. We
require some preparation. 

\begin{lem}
\label{lem:transitivity-implies-epsilon-periods}If $(X,T)$ is forward
transitive then it contains an $\varepsilon $-dense $\varepsilon $-pseudo
period.
\end{lem}
\begin{proof}
Let $x$ be a forward transitive point. Then $\{x,Tx,\ldots ,T^{m}x\}$
is $\varepsilon $-dense in $X$ for some $m$. Choose $n>m$ with
$d(T^{n}x,x)<\varepsilon $; then $x,Tx,\ldots ,T^{n-1}x$ is an $\varepsilon $-dense
$\varepsilon $-pseudo-period.
\end{proof}
The following lemma implies that approximation in $\all $ by periodic
systems, which a-priori is a property of the embedding of the system
in $Q^{\mathbb{Z}}$, is really an intrinsic property of the system.
Note that the condition in the lemma does not depend on the metric.

\begin{lem}
\label{lem:characterization-of-limits-of-cycles}Let $X\in \all _{\trans }$.
Then $X$ is the limit of finite cycles if and only if for every $\varepsilon >0$
there is an $\varepsilon $-dense $\varepsilon $-pseudo-period in
$X$. 
\end{lem}
\begin{proof}
If $\overline{x}=x_{1},\ldots ,x_{n}$ is an $\varepsilon $-dense
$\varepsilon $-pseudo-period in $X$ let $a=\theta (\overline{x})$
and let $X_{\varepsilon }\in \all $ be the cycle consisting of the
orbit of the infinite concatenation $a$ with itself. Then by \ref{cor:symbolic-approximation}
we have $d(X,X_{\varepsilon })<2\sqrt{\varepsilon }$, so $X$ is
the limit of the finite cycles systems $X_{\varepsilon }$.

Conversely suppose $X_{n}$ are finite cycles and $X_{n}\rightarrow X$.
Suppose $X_{n}$ is the orbit of the periodic point $x_{n}$ of period
$N(n)$. For each $k=0,\ldots ,N(n)-1$ let $x'_{n,k}\in X$ be one
of the points in $X$ closest to $\sigma ^{k}x_{n}$. One verifies
that $(x_{n,k})_{k=0}^{N(n)-1}$ is an $2\varepsilon $-dense $\delta (\varepsilon )$-pseudo-period
in $X$ for some $\delta (\varepsilon )\rightarrow 0$ with $\varepsilon $. 
\end{proof}
\begin{prop}
\label{prop:per-is-dense-in-trans}The finite cycles are dense in
$\all _{\trans }$
\end{prop}
\begin{proof}
Suppose $X$ is transitive and let $\varepsilon >0$. Take a recurrent
transitive point $x\in X$, which is forward recurrent by lemma \ref{lem:implications-of-forward-transitivity}.
We can therefore choose $n$ so that $x,\sigma x,\ldots ,\sigma ^{n}x$
is an $\varepsilon $-dense $\varepsilon $-pseudo-period. Apply the
previous lemma. 
\end{proof}

\subsection{\label{sub:Odometers}Odometers}

Given a factor map $f:(X,T)\rightarrow (Y,S)$, the sets $f^{-1}(y)$
for $y\in Y$ are called the \emph{fibres} (of $f$), and induce a
partition of $X$ into closed, pairwise disjoint sets. If $Y$ is
a finite cycle then the fibres are also open sets. 

The next simplest transitive systems after the finite cycles are those
systems that are determined by factors onto finite cycles. A system
$(X,T)$ is an \emph{odometer} (or \emph{adding machine}) if, for
every $\varepsilon >0$, there is a factor map $\pi $ from $X$ onto
a finite cycle $(Y,S)$ with fibres of diameter $<\varepsilon $. 

Equivalently, for every $\varepsilon >0$ there is a partition of
$X$ into closed, pairwise disjoint sets which are permuted cyclically
by $T$; these are the fibres of the factor map onto the quotient
space of the partition.

In the next few lemmas we establish some well-known properties of
odometers and define the universal odometer.

\begin{lem}
\label{lem:properties-of-odometers}Odometers are zero-dimensional
and forward transitive.
\end{lem}
\begin{proof}
The first statement is clear from the definition, since fibres of
factor maps to cycles form a closed and open basis for the topology. 

If $(X,T)$ is an odometer, we show that every point is forward transitive.
Fix $x\in X$ and let $\emptyset \neq U\subset X$ be an open and
closed set. Choose a factor onto a finite cycle with fibres of diameter
less than $\max \{d(x,y)\, :\, x\in U,y\in X\setminus U\}$. Then
each fibre is contained in $U$ or $X\setminus U$ and in particular
there are fibres contained in $U$. Since $T^{n}x$ visits this fibre
for some $n>0$, for this $n$ we have $T^{n}x\in U$. Since the closed
and open sets are a basis for the topology on $X$, the forward orbit
of $x$ is dense.
\end{proof}
It follows from lemma \ref{lem:implications-of-forward-transitivity}
that an Odometer is either a finite cycle or a Cantor system. 

\begin{lem}
\label{lem:joining-of-finite-factors}Let $(X,T)$ be a dynamical
system, let $f:X\rightarrow C_{m}$ and $g:X\rightarrow C_{n}$ be
factor maps. Then there is a factor map $h:X\rightarrow C_{k}$ such
that the maps $f,g$ factor through $h$, i.e. the fibres of $h$
refine the fibres of $f$ and of $g$, and $m,n$ divide $k$.
\end{lem}
\begin{proof}
Let $X_{0},\ldots ,X_{m-1}$and $Y_{0},\ldots ,Y_{n-1}$ be the fibres
of $f,g$ respectively. Consider the partition \[
\{X_{i}\cap Y_{j}\, :\, 0\leq i<m\, ,\, 0\leq j<n\}\]
Let $k$ be the number of these atoms. The atoms are open and closed
and $T$ acts on them by permutation; by transitivity, this permutation
is cyclic, for otherwise we could choose two atoms whose $T$-orbits
are disjoint. Let $h$ be the projection to the quotient space of
this partition with the quotient action. The relations $f=f\circ h$
and $g\circ h$ are now clear. Since $T$ is a homeomorphism and maps
fibres to fibres, the number of fibres of $h$ in each fibre of $f$
is constant, so $m|k$, and similarly $n|k$.
\end{proof}
One way to construct an odometer is as follows. Fix a sequence $(k(n))_{n\in \mathbb{N}}$
of positive integers such that $n(k)|n(k+1)$. Let $\pi _{k-1}:C_{n(k)}\rightarrow C_{n(k-1)}$
be the factor map given by $\pi _{k-1}(i)=i\bmod n(k-1)$. We obtain
a system of factor maps\[
\ldots \xrightarrow{\pi _{k+1}}C_{n(k+1)}\xrightarrow{\pi _{k}}C_{n(k)}\xrightarrow{\pi _{k-1}}C_{n(k-1)}\xrightarrow{\pi _{k-2}}\ldots \xrightarrow{\pi _{1}}C_{n(1)}\]
The inverse limit of this system is an odometer $(X,T)$ (note that
the $\pi _{k}$'s are group homomorphisms, so the inverse limit also
has a group structure; the map $T$ of $X$ is translation on this
group by the element which is the limit of $1\in C_{n(k)}$).

Let $(X^{*},T^{*})$ denote the odometer constructed as above from
the sequence $n(k)=k!$. We call this the \emph{universal odometer.}
This system may be characterized as follows:

\begin{lem}
\label{lem:universal-odometer-characterization}Up to isomorphism,
$(X^{*},T^{*})$ is the unique odometer which factors onto every finite
cycle.
\end{lem}
\begin{proof}
Fix $m\in \mathbb{N}$. By construction, $(X^{*},T^{*})$ factors
onto $C_{m!}$ and this system factors onto $C_{m}$ by reduction
modulo $m$.

Conversely suppose that $(X,T)$ is an odometer which factors onto
every finite cycle. We define a sequence $n(k)$ and factor maps $f_{k}:X\rightarrow C_{n(k)}$
by induction, as follows. Set $n(1)=1$ and $f_{1}$ maps $X$ to
a point. Suppose we are given $f_{k}:X\rightarrow C_{n(k)}$. Since
$X$ factors onto every $C_{N}$ it factors onto $C_{n(k)!}$, and
by lemma \ref{lem:joining-of-finite-factors} we can find $n(k+1)$
and a factor $f_{k+1}:X\rightarrow C_{n(k+1)}$ such that $n(k)!|n(k+1)$,
the fibres of $f_{k+1}$ are of diameter $<1/k$ and they refine the
fibres of $f_{k}$. Let $\rho _{k}$ be the factor map from $C_{n(k+1)}\rightarrow C_{n(k)}$
so that $f_{k}=\rho _{k}\circ f_{k+1}$. It now follows that $(X,T)$
is the inverse limit of the diagram\[
\ldots \xrightarrow{\rho _{k+1}}C_{n(k)}\xrightarrow{\rho _{k}}C_{n(k-1)}\xrightarrow{\rho _{k-1}}C_{n(k-2)}\xrightarrow{\rho _{k-2}}\ldots \xrightarrow{\rho _{!}}C_{n(1)}\]
Since $n(k)\, |\, n(k)!\, |\, n(k+1)$, we may interpolate $C_{n(k)!}$
between $C_{n(k+1)}$ and $C_{n(k)}$. We get \[
\ldots \rightarrow C_{n(k)!}\rightarrow C_{n(k)}\rightarrow C_{n(k-1)!}\rightarrow C_{n(k-1)}\rightarrow \ldots \rightarrow C_{n(1)}\]
It follows that $(X,T)$ is the inverse limit of \[
\ldots \rightarrow C_{n(k)!}\rightarrow C_{n(k-1)!}\rightarrow \ldots \rightarrow C_{n(1)!}\]
We may now interpolate all the other $C_{m!}$'s for $m=1,2,3\ldots $
into this sequence, and this gives us the sequence which defines $(X^{*},T^{*})$;
so $X\cong X^{*}$.
\end{proof}
Similar arguments show that every odometer can be obtained as the
inverse limit of $C_{n(k)}$'s for some sequence $n(k)$, and that
every odometer is a factor of $(X^{*},T^{*})$. This justifies the
claim of universality. We will not use these facts, and omit the proofs.

\subsection{\label{sub:genericity-of-odometer}The Universal Odometer is Generic
in $\all _{\trans }$.}

We can now bring everything together and show that the universal odometer
is generic in $\all _{\trans }$.

\begin{lem}
\label{lem:propogation-of-finite-factors}Suppose that $Y\in \all $
factors onto a finite cycle of period $n$ with fibres of diameter
$<\varepsilon $. Then the same is true of every $X\in \all $ sufficiently
close to $Y$.
\end{lem}
\begin{proof}
Suppose $Y_{0},\ldots ,Y_{n-1}$ is a partition of $Y\in \all $ into
closed sets of diameter $<\varepsilon $. Let $\delta =\min _{i\neq j}d(Y_{i},Y_{j})$.
Since the shift $\sigma $ on $Q^{\mathbb{Z}}$ is uniformly continuous
there is an $\eta $ so that if $d(x',x'')\leq \eta $ then $d(\sigma x',\sigma x'')\leq \delta /3$.
We may assume that $\eta <\delta $ and $2\eta +\diam Y_{i}<\varepsilon $
for $i=0,\ldots ,n-1$. Let $X\in \all $ be a system with $d(X,Y)<\eta $.
Let \[
X_{i}=\{x\in X\, :\, d(x,Y_{i})<\eta \}\]
this is a closed partition of $X$ into sets of diameter $<2\eta +\max _{i}\diam Y_{i}<\varepsilon $.
For every point $x\in X_{i}$ there is a $y\in Y_{i}$ with $d(x,y)<\eta $,
so $d(\sigma x,\sigma y)<\delta /3$, and since $\sigma y\in Y_{i+1}$
and $\eta <\delta /3$ we conclude that $d(\sigma x,Y_{j})>\eta $
for $j\neq i$, so $\sigma x\in X_{i+1}$. Thus $\sigma $ permutes
the $X_{i}$ cyclically. This completes the proof. 
\end{proof}
\begin{lem}
\label{lem:density-of-hi-period-cycles}For every $n$, the finite
cycles of period divisible by $n$ are dense in $\all _{\trans }$.
\end{lem}
\begin{proof}
Fix $n\in \mathbb{N}$ and let $Y\in \all $ be a finite cycle of
period $k$. Then $Y$ consists of the orbit of a point $y\in Q^{\mathbb{Z}}$
of period $k$, that is, $y(i+k)=y(i)$ for all $i\in \mathbb{Z}$
and $k$ is the least positive integer with this property. Since $Q$
has no isolated points we an perturb the coordinates of $y$ by a
small amount in a manner which has period $kn$. We obtain a point
$y'$, which can be made arbitrarily close to $y$ in $d_{\infty }$,
with period $kn$. The orbit of this point is a finite cycle $Y'$
of period $kn$, and $Y'\rightarrow Y$ as $y'\rightarrow y$. 

We have shown that the cycles of period divisible by $n$ are dense
among the finite cycles, and since by \ref{prop:per-is-dense-in-trans}
the latter are dense in $\all $, we are done.
\end{proof}
\begin{thm}
\label{thm:genericity-of-odometer}The isomorphism class of $(X^{*},T^{*})$
is a dense $G_{\delta }$ set in $\all _{\trans }$.
\end{thm}
\begin{proof}
Let $\mathcal{U}_{n}\subseteq \all $ be the set of transitive systems
which factor onto a finite cycle of period divisible by $n$ with
fibres of diameter $<1/n$. By lemma \ref{lem:propogation-of-finite-factors}
this set is open in $\all $, and hence in $\all _{\trans }$, and
by lemma \ref{lem:density-of-hi-period-cycles} this set is dense.
By lemma \ref{lem:universal-odometer-characterization}, $\cap _{n=1}^{\infty }\mathcal{U}_{n}$
is the isomorphism class of $(X^{*},T^{*})$.
\end{proof}
In particular, this establishes that there is a dense set of Cantor
systems in $\all _{\trans }$. Together with lemma \ref{pro:trans-is-Gdelta-in-G}
and \ref{thm:trans-is-Gdelta-in-S} we have verified the hypotheses
of the correspondence theorem, proving theorem \ref{thm:transitive-systems}.

\section{\label{sec:total-transitivity}The Space of Totally Transitive Systems}

As we have seen, among the transitive systems periodicity phenomena
are quite dominant. These are absent from the class of totally transitive
systems. This class, which contains a rich variety of dynamics and
a few surprises, will occupy us for the remainder of this paper.

\subsection{\label{sub:Total-Transitivity}Total Transitivity }

A system $(X,T)$ is \emph{totally transitive} if $T^{n}$ is transitive
for every $n>0$. The class of totally transitive systems is denoted
by $\ttrans $. 

We remark that in this definition it is not important whether we interpret
transitivity as forward- or bi-transitivity; the result is the same.
Indeed, if $T^{n}$ is forward transitive for each $n$, then it is
bi-transitive for every $n$. conversely, suppose $(X,T)$ is a system
with $T^{n}$ bi-transitive for each $n$ and suppose $x\in X$ is
an isolated point. Since $\{x\}$ is open it must be a bi-transitive
point for each $T^{n}$. In particular there is a $k\in \mathbb{Z}$
with $T^{2k}x=Tx$. But this is impossible because $T$ is injective.
Therefore $X$ was perfect, so bi-transitivity of $T^{n}$ implies
forward transitivity of $T^{n}$ by lemma \ref{lem:transitivity-in-perfect-spaces}.

Let us say a few words about the relation between $\all _{\trans }$
and $\all _{\ttrans }$. The totally transitive systems are not dense
in $\all _{\trans }$ since the set of systems in $\all _{\trans }$
with periodic factors is open (lemma \ref{lem:propogation-of-finite-factors})
and dense (in particular, nonempty), and a totally transitive system
cannot have a periodic factor. Note however that totally transitivity
is not equivalent to the nonexistence of periodic factors

The closure $\overline{\all }_{\ttrans }$ is not contained in $\all _{\trans }$.
For instance let $0,1\in Q$ be distinct points and let $X_{n}$ be
the system over the alphabet $\{0,1\}$ in which every maximal block
of consecutive $0$'s or consecutive $1$'s has length at least $n$.
Then $X_{n}$ are seen to be totally transitive, but $\lim X_{n}$
exists and is the simple system generated by the very simple points
points $(\ldots ,0,0,1,1,\ldots )$ and $(\ldots ,1,1,0,0,\ldots )$.
This system is not transitive.

$\all _{\ttrans }$ is not even relatively closed in $\all _{\trans }$.
Consider for example the systems $X_{n}\subseteq \{0,1\}^{\mathbb{Z}}$
which contain all sequences of $0$'s and $1$'s such that if two
occurrences of $1$'s occur at distance $k\leq n$ then $k$ is even.
These systems are totally transitive, and $X_{n}\rightarrow X\subseteq \{0,1\}^{\mathbb{Z}}$,
where $X$ contains all sequences of $0$'s and $1$'s in which the
distance between every two $1$'s is even. $X$ is not totally transitive.

The following is needed in order to apply the correspondence theorem
to $\ttrans $:

\begin{prop}
\label{prop:ttrans-is-Gdelta}$\all _{\ttrans }$ is a $G_{\delta }$
in $\all $ and $\mathcal{H}_{\ttrans }$ is a $G_{\delta }$in $\mathcal{H}$.
\end{prop}
\begin{proof}
For each $r$ one can imitate the proof that $\all _{\trans }$ is
a $G_{\delta }$ and obtain that the space \[
\mathcal{T}_{r}=\{X\in \all \, :\, (X,\sigma ^{r})\textrm{ is transitive}\}\]
is a $G_{\delta }$. By definition, $\all _{\ttrans }=\cap _{r}\mathcal{T}_{r}$.
The proof in $\mathcal{H}$ is similar.
\end{proof}
We conclude with an important combinatorial property of totally transitive
systems, upon which we will rely heavily in the sequel.

\begin{prop}
\label{pro:ttrans-implies-independent-epsilon-periods}Let $(X,T)\in \ttrans $.
Then for every $\varepsilon >0$ there exist two relatively prime
numbers $n,m$ and two $\varepsilon $-dense $\varepsilon $-pseudo
periods $x_{1},\ldots ,x_{n},$ and $x'_{1},\ldots ,x'_{m}$ in $X$
with the same starting point, $x_{1}=x'_{1}$.
\end{prop}
\begin{proof}
Suppose $(X,T)$ is totally transitive, let $\varepsilon >0$, and
consider the set $I$ of integers $n$ such that there exists an $\varepsilon $-dense
$\varepsilon $-pseudo period $x_{1},x_{2},\ldots ,x_{n}$ in $X$.
Let $r$ be the least common divisor of $I$. If $r=1$ we are done.
Otherwise $r\neq 1$ and $r|n$ for every $n\in I$. Fix a transitive
point $x_{1}\in X$ and for $i=0,\ldots ,r-1$ let $A_{i}\subseteq X$
consist of the points $y$ such that there is an $\varepsilon $-pseudo
orbit $x_{1},x_{2},\ldots ,x_{n}$ with $n\equiv i(\mod r)$ and $x_{n}=y$. 

The $A_{i}$ are clearly open, and their union is all of $X$ because
$x_{1}$ is a transitive point.

The $A_{i}$ are also disjoint: if $x'\in A_{i}\cap A_{j}$ and $i<j$
then we can construct an $\varepsilon $-pseudo period starting at
$x_{1}$ of length $j-i(\mod r)$ by choosing a $\varepsilon $-pseudo-orbits
$x_{1},x'_{2},\ldots ,x'_{m}=x',x'_{m+1},\ldots ,x'_{M}=x_{1}$ with
$m=i(\mod r)$ and $x'_{1},\ldots ,x'_{n}=x',x'_{n+1},\ldots ,x'_{N}$
with $n=j(\mod r)$, and forming the $\varepsilon $-pseudo-period
$x_{1},x'_{2},\ldots ,x'_{m},x''_{n+1},\ldots ,x''_{N}$. This contradicts
the definition of $r$.

Since $TA_{i}=A_{i+1(\mod r)}$ the projection taking $x\in X$ to
the unique $A_{i}$ to which it belongs is a factor map from $X$
to the periodic system $(\{A_{0},\ldots ,A_{r-1}\},T)$. This contradicts
the total transitivity of $X$.
\end{proof}

\subsection{\label{sub:ttrans-rokhlin-property}The Rohlin Property }

There is a classical theorem of Rohlin in ergodic theory, stating
that the isomorphism class of every aperiodic ergodic system is dense
in the automorphism group of a Lebesgue space (aperiodicity means
that the set of periodic points has measure zero). We next discuss
the topological analogue of this, i.e. when the isomorphism class
of a topological system $(X,T)$ is dense in $\all _{\ttrans }$. 

There are two obvious obstructions to this. One is the existence of
periodic points. Suppose that $(X,T)$ has a periodic point $x$ of
period $k$ and that $X_{n}\rightarrow Y$ with $(X_{n},\sigma )\cong (X,T)$.
Let $x_{n}\in X_{n}$ have period $k$. If $y$ is any accumulation
point of $x_{n}$ then $y\in Y$, and $\sigma ^{k}y=y$. In particular,
the closure of systems with periodic points cannot contain systems
without periodic points, such as infinite minimal systems.

Another obstruction is topological. Clearly, if $X$ is connected
and $Y$ disconnected, then one cannot approximate $Y$ in the Hausdorff
metric by homeomorphic images of $X$.

It turns out that these are essentially the only obstructions. We
will show that the isomorphism class of any zero dimensional system
without periodic points is dense in $\all _{\ttrans }$. In the proof
we use the following lemma, which is a variation on a lemma of Krieger
\cite{K}.

\begin{lem}
Let $(X,T)$ be a zero-dimensional system without periodic points.
Then for every $k$ there is a closed and open subset $U\subseteq X$
such that $U,TU,\ldots ,T^{k}U$ are pairwise disjoint, and $X=\cup _{i=-k}^{k}T^{i}U$.
\end{lem}
\begin{proof}
Since no point in $X$ is periodic, for $x\in X$ the points $x,Tx,\ldots ,T^{k}x$
are distinct, so we may choose a closed and open set $E_{x}$ containing
$x$ and $E_{x},TE_{x},\ldots ,T^{k}E_{x}$ are pairwise disjoint.
The collection $\{E_{x}\}_{x\in X}$ covers $X$; choose a finite
subcover $E_{1},\ldots ,E_{m}$. Set $E'_{1}=E_{1}$ and \[
E'_{i+1}=E'_{i}\cup (E_{i+1}\setminus \cup _{j=-k}^{k}E'_{i})\]
One now checks that $U=E'_{m}$ satisfies our requirements.
\end{proof}
\begin{thm}
\label{thm:dense-isomorphism-class}Let $(X,T)$ be zero-dimensional
system without periodic points. Then the set of systems isomorphic
to $X$ is dense in the totally transitive systems, i.e. $\all _{\ttrans }\subseteq \overline{\all }_{[(X,T)]}$
in $\all $. In particular, if $(X,T)$ is totally transitive then
its isomorphism class is dense in $\all _{\ttrans }$.
\end{thm}
\begin{proof}
Let $Y\in \all _{\ttrans }$ be totally transitive and $\varepsilon >0$;
we will find a system $X'\in \all $ isomorphic to $(X,T)$ and with
$d(X',Y)<\varepsilon $.

The proof has two steps. First, use the lemma to get a symbolic factor
of $X$ which is within $\varepsilon $ of $Y$, and then modify this
map to make it an embedding of $X$ into $Q^{\mathbb{Z}}$ without
moving the image more than $\varepsilon $.

Here is the proof of the first step. Select two $\varepsilon $-dense
$\varepsilon $-pseudo-periods in $Y$ with relatively prime lengths
$m,n$ and common starting point: $\overline{y}'=y'_{1},\ldots ,y'_{m}$
and $\overline{y}''=y''_{1},\ldots ,y''_{n}$. Let $a'=\theta (\overline{y}')$
and $a''=\theta (y'')$ (recall the notation of section \ref{sec:symbolic-approximation}). 

There exists an integer $k_{0}$ such that every integer $k>k_{0}$
can be written as $k=rm+sn$ for integers $r,s\in \mathbb{N}$. For
each $k>k_{0}$ we can therefore fix a word $a_{k}$ of length $k$
consisting of concatenations of $a',a''$.

Using the lemma, select a closed and open set $U\subseteq X$ such
that if $x\in U$ then $T^{i}u\notin U$ for $i=1,\ldots ,k_{0}$,
but $\cup _{i=-k_{0}}^{k_{0}}T^{i}U=X$; hence for every $x\in X$
we have $T^{i}x\in U$ for some $0\leq i\leq k_{0}$ and $T^{j}x\in U$
for some $-k_{0}\leq j<0$.

We now define a continuous map $f:X\rightarrow Q^{Z}$ using visits
to $U$ as {}``markers''. For $x\in X$, let $I\subseteq \mathbb{Z}$
be the set of times $i\in \mathbb{Z}$ such that $T^{i}x\in U$. By
choice of $U$ the set $I$ is bounded neither above nor below, and
the gap between consecutive times in $I$ is at most $M=2k_{0}$.
By choice of $U$ we also know that if $i,j\in I$ and $i\neq j$
then $|i-j|>k_{0}$. Fix an ordering $I=\{\ldots <i(-1)<i(0)<i(1)<i(2)<\ldots \}$
and form the word $f(x)\in Q^{\mathbb{Z}}$ such that at index $i(m)$
we see the word $a_{i(m+1)-i(m)}$. The map $x\mapsto f(x)$ obviously
satisfies $f(Tx)=\sigma f(x)$, and $f$ is continuous since the symbol
$f(x)(0)$ is determined by the values of the minimal $i,j\in \mathbb{N}$
such that $T^{i}x\in U$ and $T^{-j}x\in U$; and these vary continuously
in $x$ because $U$ is open and closed. We see that $f$ is a factor
map onto its image. The image $f(X)$ is constructed from the language
$\{a_{k}\}_{k\geq k_{0}}$, so it is constructed from the language
$\{a',a''\}$, which by the choice of $a',a''$ and corollary \ref{cor:symbolic-approximation}
gives gives $d(f(X),Y)\leq 2\sqrt{\varepsilon }$ .

In general $f(X)\not \equiv X$ because $f$ need not be injective.
The second step is to perturb $f$ to make it so. Let\[
Q_{0}=\{a'(i),a''(j)\, :\, i=1,\ldots ,m\, ,\, j=1,\ldots ,n\}\subseteq Q\]
 be the list of symbols appearing in $a',a''$ For each $z\in Q_{0}$
let $X_{z}\subseteq Q$ be a homeomorphic copy of $X$ such that $d(z,x)<\varepsilon $
for every $x\in X_{z}$, and such that the $X_{z}$'s are pairwise
disjoint. This can be done because $Q_{0}$ is finite and $Q$ is
perfect. For each $z\in Q_{0}$ let $g_{z}:X\rightarrow X_{z}$ be
a fixed homeomorphism.

Define a map $h:X\rightarrow Q^{\mathbb{Z}}$ by\[
h(x)(i)=g_{f(x)(i)}(\sigma ^{i}x)\]
$h$ is clearly continuous, and commutes with the shift since \[
h(\sigma x)(i)=g_{f(\sigma x)(i)}(\sigma ^{i}(\sigma x))=g_{f(x)(i+1)}(\sigma ^{i+i}x)=h(x)(i+1)\]
 The map $h$ is is an injection, since if $y\in h(X)$ then $y(0)\in X_{z}$
for some $z\in Q_{0}$, implying that $g_{z}^{-1}(y)$ is the unique
preimage of $y$ under $h$. Hence $h$ is an isomorphism onto its
image.

Finally, since $d(z,g_{z}(x))<\varepsilon $ for every $x\in X$ and
$z\in Q_{0}$, we see that $d_{\infty }(f(x),h(x))<\varepsilon $.
Therefore $d(f(X),h(X))<\varepsilon $, so\[
d(h(X),Y)\leq d(h(X),f(X))+d(f(X),Y)\leq \varepsilon +2\sqrt{\varepsilon }\]
So $h(X)$ is the desired system.
\end{proof}
\begin{cor}
The correspondence theorem and the zero-one law hold in $\ttrans $.
\end{cor}
\begin{proof}
There exist totally transitive Cantor systems without periodic points;
any minimal totally transitive Cantor system, for instance. We conclude
that in $\all _{\ttrans }$ there is a dense set of Cantor systems.
Together with proposition \ref{prop:ttrans-is-Gdelta}, this completes
the hypotheses of the correspondence theorem for $\all _{\ttrans }$and
$\mathcal{H}_{\ttrans }$.
\end{proof}
\begin{cor}
The isomorphism class of any Cantor system without periodic points
is dense in $\mathcal{H}_{\ttrans }$. In particular the zero-one
law (theorem \ref{thm:zero-one-law}) holds.
\end{cor}
\begin{proof}
Let $(X,T)$ be a minimal totally transitive Cantor system. Set $\mathbf{Q}=[(X,T)]$
and apply part (2) of the correspondence theorem to conclude that
$\mathbf{Q}$ is dense in $\mathcal{H}_{\ttrans }$. The zero one
law follows from theorem \ref{thm:zero-one-law}.
\end{proof}

\subsection{Minimality, Unique Ergodicity and Entropy}

A dynamical system $(X,T)$ is minimal if it has no nontrivial subsystems,
or equivalently, if every point is (forward) transitive. Minimality
can also be characterized by the property that for every nonempty
open set $U\subseteq X$ there is an $N$ such that for every $x\in X$
we have $T^{i}x\in U$ for some $i\in \{1,\ldots ,N\}$. For a proof
of this, see \cite[Theorem 5.1]{Wal81}. The class of minimal systems
is denoted $\mn $. 

A systems $(X,T)$ is \emph{uniquely ergodic} if there is a unique
$T$-invariant Borel probability measure on $X$. The class of uniquely
ergodic systems is denoted $\ue $. Unique ergodicity can be characterized
by the property that for every continuous function $f\in C(X)$ and
every $\varepsilon >0$ there is an $N$ such that for every $x,y\in X$,\[
|\frac{1}{N}\sum _{i=1}^{N}f(T^{i}x)-\frac{1}{N}\sum _{i=1}^{N}f(T^{i}y)|<\varepsilon \]
This follows from \cite[Theorem 6.19]{Wal81}.

\begin{prop}
$\mn $ and $\ue $ are dense $G_{\delta }$ subsets of $\all _{\ttrans }$.
\end{prop}
\begin{proof}
Density of both classes follows from theorem \ref{thm:dense-isomorphism-class}. 

Let $\{U_{i}\}_{i=1}^{\infty }$ be a countable basis for the topology
of $\all $. As in the proof of \ref{thm:trans-is-Gdelta-in-S}, for
each $i$ and $N$ the set of $X\in \all $ such that\[
\mathcal{U}_{i,N}=\{X\in \all \, :\, X\cap U_{i}\neq \emptyset \textrm{ and }X\subseteq \bigcup _{k=0}^{N}\sigma ^{k}U_{i}\}\]
is open in $\all $, and so is the union $\mathcal{U}_{i,N}$ over
$N$. Since the set \[
\mathcal{V}_{i}=\{X\in \all \, :\, X\cap U_{i}=\emptyset \}\]
is closed it is a $G_{\delta }$. Therefore \[
\mathcal{W}_{i}=\mathcal{V}_{i}\cup \bigcup _{N=1}^{\infty }\mathcal{U}_{i,N}=\{X\in \all \, :\, X\cap U_{i}\neq \emptyset \Rightarrow X\subseteq \bigcup _{k=0}^{N}\sigma ^{k}U_{i}\}\]
 is a $G_{\delta }$. Now $\all _{\mn }=\cap _{i=1}^{\infty }\mathcal{W}_{i}$
is a $G_{\delta }$.

A similar argument shows that $\all _{\ue }$ is a $G_{\delta }$.
Briefly, for each $i$ and $n,N$ one shows that the set\[
\mathcal{U}_{i,n,N}=\{X\in \all \, :\, |\frac{1}{N}\sum _{i=1}^{N}f_{i}(\sigma ^{i}x)-\frac{1}{N}\sum _{i=1}^{N}f_{i}(\sigma ^{i}y)|<\frac{1}{n}\textrm{ for all }x,y\in X\}\]
 is open. Then $\all _{\ue }=\cap _{i,n}\cup _{N}\mathcal{U}_{i,n,N}$
is a $G_{\delta }$.

The proofs in $\mathcal{H}$ are similar and slightly easier, so we
omit them.
\end{proof}
Our next result requires familiarity with Bowen's definition of entropy
(see \cite{Wal81}). 

\begin{prop}
The class of zero-entropy systems is a dense $G_{\delta }$ in $\all _{\ttrans }$. 
\end{prop}
\begin{proof}
Density is clear from theorem \ref{thm:dense-isomorphism-class},
since there exists minimal totally transitive Cantor systems with
zero entropy.

To see that the class of zero-entropy systems is a $G_{\delta }$,
let $\mathcal{U}_{n}\subseteq \all $ be the set of systems $X$ such
that for some $k$ the number of $1/n$-separated words of length
$k>n$ in $X$ is less than $2^{k/n}$. One verifies that the $\mathcal{U}_{n}$
are open, and their intersection consists exactly of the zero entropy
systems.
\end{proof}

\subsection{Connections with the Space of Invariant Measures}

In this section we will prove a partial correspondence theorem relating
the genericity of a dynamical property $\mathbf{P}$ in the measure-theoretic
category to the genericity in the class of totally transitive topological
systems which support an invariant measure in $\mathbf{P}$.

Recall that $\mathcal{M}$ is the space of shift-invariant Borel probability
measures on $Q^{\mathbb{Z}}$ with the weak-{*} topology. 

\begin{thm}
\label{thm:trans-category-correspondence}Let $\mathbf{P}$ be a dynamical
property in the measure theoretic category and suppose that $\mathcal{M}_{\mathbf{P}}$
is a dense $G_{\delta }$ in $\mathcal{M}$. Let $\widetilde{\mathbf{P}}$
be the class of topological dynamical systems which support a measure
from $\mathbf{P}$. Then $\widetilde{\mathbf{P}}$ is generic in $\all _{\ttrans }$.
\end{thm}
\begin{proof}
For a system $X\in \mn \cap \ue $, denote the unique invariant measure
by $\mu _{X}$. Let $m:\all _{\mn \cap \ue }\rightarrow \mathcal{M}$
be the map $X\mapsto \mu _{X}$. One may verify that this map is continuous
(but note that its image is meagre in $\mathcal{M}$, since generically
a measure in $\mathcal{M}$ has global support). Thus $m^{-1}(\mathcal{M}_{\mathbf{P}})\subseteq \all _{\mn \cap \ue }$
is a $G_{\delta }$ in $\mn \cap \ue $ and hence in $\ttrans $;
and $m^{-1}(\mathcal{M}_{\mathbf{P}})\subseteq \all _{\widetilde{\mathbf{P}}}$. 

Since zero entropy and weak mixing (in the ergodic sense) are generic
in $\mathcal{M}$ and $\mathcal{M}_{\mathbf{P}}$ is generic in $\mathcal{M}$,
there is a weak mixing zero entropy system $(Y,S,\nu )\in \mathcal{M}_{\mathbf{P}}$.
By Krieger's generator theorem, there is a minimal uniquely ergodic
symbolic system on two symbols $(Z,T)$ such that $(Z,T,\mu _{Z})\cong (Y,S,\nu )$
in the measure-theoretic category. Thus there is a system $X\in \all _{\mn \cap \ue }$
with $(X,\sigma )\cong (Z,T)$ in the topological category. Since
$(X,\sigma ,\mu _{X})$ is measure-theoretically weak mixing, $(X,\sigma )$
is totally transitive. Applying theorem \ref{thm:dense-isomorphism-class}
we see that the isomorphism class of $(X,\sigma )$ is dense in $\all _{\mn \cap \ue \cap \ttrans }$.
But $X\in m^{-1}(\mathcal{M}_{\mathbf{P}})$, and so is every system
isomorphic to $X$. Thus $m^{-1}(\mathcal{M}_{\mathbf{P}})$ is dense
in $\mn \cap \ue \cap \ttrans $, so $\widetilde{\mathbf{P}}$ is
generic there.
\end{proof}
\begin{cor}
The systems supporting a global weak mixing rigid invariant measure
are generic in $\ttrans $, and in particular systems supporting a
strong mixing measure are exotic there.
\end{cor}
We note that the relation between \emph{topological} weak and strong
mixing is somewhat different; see section \ref{sec:mixing} below.

For the proof of theorem \ref{thm:trans-category-correspondence}
it was necessary to assume that $\mathbf{P}$ is a dense $G_{\delta }$
in $\mathcal{M}$. It would be nice to weaken this assumption to the
weaker one that $\mathbf{P}$ is generic in $\mathcal{M}$. The problem
is that in this case $\mathcal{M}_{\mathbf{P}}$ does indeed contain
a dense $G_{\delta }$ subset $\mathcal{M}_{\mathbf{P}}^{*}\subseteq \mathcal{M}_{\mathbf{P}}$,
but this $G_{\delta }$ set may not be saturated with respect to the
isomorphism relation, and in the above proof we cannot conclude that
the system $(X,\sigma )$ is in $m^{-1}(\mathcal{M}_{\mathbf{P}}^{*})$.
We do not know if the theorem holds under weaker hypotheses.

\section{\label{sec:Disjointness}Disjointness}

For dynamical systems $(X,T),(Y,S)$ let $T\times S$ denote the homeomorphism
of $X\times Y$ given by $(T\times S)(x,y)=(Tx,Sy)$. 

Let $(X,T),(Y,S)$ be bi-transitive dynamical systems. A \emph{joining}
of $(X,T),(Y,S)$ is a bi-transitive subsystem $Z\subseteq X\times Y$
whose projection to the first coordinate is onto $X$ and to the second
coordinate is onto $Y$. Systems $(X,T),(Y,S)$ are \emph{disjoint}
if their only joining is the product system \emph{$(X\times Y,T\times S)$.}

The notion of disjointness, along with an analogous notion in the
measure-preserving category, was introduced by Furstenberg in \cite{F}
as a tool in the classification of dynamical systems and has proved
a very successful one. In \cite{dJ}, del Junco showed that for any
automorphism $T$ of a Lebesgue space the set of automorphisms measure-theoretically
disjoint from $T$ is residual in the coarse topology on the space
of automorphisms. We prove the following analogue of this:

\begin{thm}
\label{thm:disjoint-is-residual}Let $X\in \ttrans $. Then a generic
system in $\all _{\ttrans }$ is disjoint from $X$.
\end{thm}
For the proof we will need a few simple facts, which  we provide for
completeness:

\begin{lem}
In order for bi-transitive systems $(X,T),(Y,S)$ to be disjoint it
suffices that for every two bi-transitive points $x\in X$ and $y\in Y$,
the point $(x,y)$ is a bi-transitive point for $X\times Y$. 
\end{lem}
\begin{proof}
Suppose this holds and $Z\subseteq X\times Y$ is a joining; then
it has a bi-transitive point $(x,y)\in Z$. Since the closure of the
full orbit of $x$ is the projection of the closure of the full orbit
of $(x,y)$ in $Z$, and $Z$ is a joining, $x$ is a bi-transitive
point for $X$. Similarly, $y$ is bi-transitive for $Y$, and by
our assumption, $Z=X\times Y$.
\end{proof}
\begin{lem}
If $(X,T)$ is totally transitive and if $x\in X$ is a bi-transitive
point, then $x$ is bi-transitive for $(X,T^{n})$.
\end{lem}
\begin{proof}
Fix $n$ and let $X_{0}$ be the closure of $\{T^{kn}x\}_{k\in \mathbb{Z}}$.
Then $Y=X_{0}\cup TX_{0}\cup \ldots \cup T^{n-1}X_{0}$ is closed
and contains the two-sided orbit of $x$ so $Y=X$. By Baire's theorem
one of the translates $T^{i}X_{0}$ must have nonempty interior so
this is true of $X_{0}$. Since $T^{n}X_{0}=X_{0}$ and in particular
$T^{n}$ preserves the interior of $X_{0}$ we conclude from the transitivity
of $(X,T^{n})$ that the interior of $X_{0}$ is dense in $X$ so
$X_{0}=X$; as desired.
\end{proof}
\begin{lem}
Every totally transitive system is disjoint from every finite cycle
. 
\end{lem}
\begin{proof}
Let $(X,T)$ be totally transitive, and $(Y,S)$ a finite cycle with
period $k$. Let $x\in X$ be a transitive point and fix $y\in Y$.
We have $(T\times S)^{k}(x,y)=(T^{k}(x),y)$, so that the orbit closure
of $(x,y)$ under the map $(T\times S)^{k}$ is $X\times \{y\}$ But
then the orbit closure of $(x,y)$ under $T\times S$ is the union
$\cup _{i=0}^{k-1}X\times \{S^{i}y\}=X\times Y$. This implies disjointness. 
\end{proof}
A similar argument shows that every odometer is disjoint from every
totally transitive system.

\begin{proof}
(of theorem \ref{thm:disjoint-is-residual}) Let us say that $Y\in \all _{\ttrans }$
is $\varepsilon $-disjoint from $X$ if, for every bi-transitive
point $x\in X$ and every bi-transitive point $y\in Y$, the full
orbit closure of $(x,y)$ in $X\times Y$ is $\varepsilon $-dense
in $X\times Y$ with respect to the metric \[
d((x',y'),(x'',y''))=\max \{d(x',y'),d(x'',y'')\}\]
 on $X\times Y$. If a system $Y$ is $\varepsilon $-disjoint from
$X$ for every $\varepsilon >0$ then it is disjoint from $X$.

Let $U_{r}\subseteq \all _{\ttrans }$ be the set of systems $1/r$-disjoint
from $X$; we will complete the proof by showing that $U_{r}$ contains
an open dense set.

Fix $r$ and $Y\in \all _{\ttrans }$. It suffices to show that for
every $\varepsilon >0$ there exists a system $Z\in \all _{\ttrans }$
with $d(Y,Z)<\varepsilon $ such that some neighborhood of $Z$ in
$\all _{\ttrans }$ is contained in $U_{r}$. So fix $\varepsilon >0$
and a positive $\delta $ to be determined later. Choose two $\delta $-dense
$\delta $-pseudo periods $\overline{a}=a_{1},\ldots ,a_{m}$ and
$\overline{b}=b_{1},\ldots ,b_{n}$ in $Y$ with $m,n$ relatively
prime and $a_{1}=b_{1}$. Let $a=\theta (\overline{a})\, ,\, b=\theta (\overline{b})$,
and assume that $\{a,b\}$ has the unique parsing property; see remark
after definition \ref{def:unique-parsing}.

Fix a (forward) transitive point $x_{0}\in X$ and let $A$ be the
periodic system whose single point up to translation is the bi-infinite
concatenation of $a$'s; denote this point by $a^{*}$. Now $X,A$
are disjoint, because $X$ is totally transitive and $A$ periodic,
so the orbit of $(x_{0},\sigma ^{i}a^{*})$ is dense in $X\times A$
for $i=1,\ldots ,m$. Therefore there is an integer $k_{a}$ such
that for each $i=1,\ldots ,m$, the first $k_{a}$ points on the orbit
of $(x,\sigma ^{i}a^{*})$ are $\delta $-dense in $X\times A$. Let
$b^{*},B$ be defined in the same way; there is an integer $k_{b}$
such that for every $j=1,\ldots ,n$ the first $k_{b}$ points on
the orbit of $(x_{0},\sigma ^{j}b^{*})$ are $\delta $-dense in $X\times B$.
Let $k'=\max \{k_{a},k_{b}\}$. 

Set $x_{1}=T^{k'}x$; what we have so far is that for $i=1,\ldots ,m$
the first $k'$ points in the \emph{backward} orbit of $(x_{0},\sigma ^{i}a^{*})$
is $\delta $-dense in $X\times A$ and similarly for $j=1,\ldots ,n$
the first $k'$ point on the backward orbit of $(x_{0},\sigma ^{j}b^{*})$
are $\delta $-dense in $X\times B$. Now $x_{1}$ is still a (forward)
transitive point for $X$ so by the same reasoning as in the previous
paragraph there is a $k''$ such that the first $k''$ points on the
orbits of $(x_{1},T^{i}a^{*}),(x_{1},\sigma ^{j}b^{*})$ are $\delta $-dense
in $X\times A,X\times B$ respectively, for every $i=1,\ldots ,m$
and $j=1,\ldots ,n$. 

Let $k=\max \{k',k''\}$. The point $x_{1}$ has the property that
for any $u\in A$ and any $v\in B$ the first $k$ iterates of $(x_{1},u)$
and $(x_{1},v)$, in either direction, are $\delta $-dense in $X\times Y$.
This remains true for any $x$ close enough to $x_{1}$; let $V\subseteq X$
be a small neighborhood of $x_{1}$ so that every $x\in V$ has this
property. 

Choose two integers $M,N>2k$ such that $mM,nN$ are relatively prime
and let $Z\in \ttrans $ be any system constructed from $\{a^{M},b^{N}\}$,
where $a^{M}$ is the $M$-fold concatenation of $a$ and $b^{N}$
the $N$-fold concatenation of $b$. By corollary \ref{cor:symbolic-approximation},
for $\delta $ small enough we will have $d(Y,Z)<\varepsilon $; we
claim that in addition $Z$ is $1/2r$-disjoint from $X$.

To see this, let $x\in X$ and $z\in Z$ be a bi-transitive points.
We must show that the full orbit of $(x,z)$ is $1/2r$-dense in $X\times Z$.
For some $i$, $T^{i}x\in V$. Since it suffices to show that the
full orbit of $(T^{i}x,\sigma ^{i}z)$ is $1/2r$-dense in $X\times Z$,
we may assume that $x\in V$ to begin with.

Parse $z$ into words of type $a^{M}$ or $b^{N}$. The $0$-th coordinate
of $z$ is in an $a^{M}$ block or in a $b^{N}$ block. In the first
case since $M\geq 2k$ either the $k$-block starting at $0$ in $z$,
or the $k$-block ending at $0$ in $z$, looks like a concatenation
of $a$'s. It follows that, for small enough choice of $\delta $,
the first $k$ points of either the forward or backward orbit of $(x,z)$
is $1/2r$-dense in $X\times Z$ (we assume that $\delta $ was chosen
small enough to overcome any edge effects arising from the part of
$z$ outside the $a^{M}$-block). The same argument applies in case
the $0$-th coordinate lies in a $b^{N}$ block; thus $Z$ is $1/2r$-disjoint
from $X$.

We will complete the proof by showing that any $W\in \all _{\ttrans }$
close enough to $Z$ in $\all $ is $1/r$-disjoint from $X$. Let
$W\in \all _{\ttrans }$ with $d(Z,W)<\delta $ and let $x\in X$
and $w\in W$ be bi-transitive points. 

The system $W'=\tau _{Z}(W)$ is a totally transitive system (since
it is a factor of the totally transitive system $W$) and assuming
$\delta $ is small enough $W'$ is constructed from $\{a^{M},b^{N}\}$
(see the discussion at the end of section \ref{sec:symbolic-approximation}).
Thus the argument carried out for $Z$ applies to $W'$ as well, so
$W'$ is $1/2r$-disjoint from $X$. Since $w'=\tau _{Z}(w)$ is a
bi-transitive point in $W'$, the full orbit of $(x,w')$ in $X\times W'$
is $1/2r$-dense in $X\times W'$. 

Recalling the properties of the projection $\tau _{Z}$ from section
\ref{sec:symbolic-approximation}, we have that $d_{\infty }(w,\tau _{Z}(w))<\delta $
and thus\[
d((T\times \sigma )^{i}(x,w),(T\times \sigma )^{i}(x,\tau _{Z}(w)))<\delta \]
 for every $i\in \mathbb{Z}$, so we conclude that the full orbit
of $(x,w)$ is $(1/2r+\delta )$-dense in $X\times W'$. Using the
fact that \[
d(W,\tau _{Z}(W))<\delta \]
we see that the full orbit of $(x,w)$ is $(1/2r+2\delta )$-dense
in $X\times W$. Thus $W$ is $1/r$-disjoint from $X$, and we are
done.
\end{proof}
\begin{thm}
If $\mathcal{C}\subseteq \all _{\ttrans }$ is compact (or $\sigma $-compact)
then a generic system in $\all _{\ttrans }$ is disjoint from every
$X\in \mathcal{C}$. 
\end{thm}
The proof is essentially the same as before. One must choose $M,N$
in such a way that for some $k$ and for every $X\in \mathcal{C}$
there is a transitive point $x_{X}\in X$ and a neighborhood $V_{X}\subseteq X$
of $x_{X}$ such that for every $x\in V_{X}$ the first $k$ iterates
of $(x,a^{*})$ in $X\times A$ is $\delta $-dense in $X\times A$,
and similarly for $(x,b^{*})$. This can be done because $\mathcal{C}$
is compact..

\section{\label{sec:mixing}Mixing}

\subsection{Weak Mixing}

A dynamical system $(X,T)$ is \emph{weakly mixing} if $(X\times X,T\times T)$
is transitive. The class of weakly mixing systems is denoted by $\wm $. 

In the definition of weak mixing we can replace transitivity of $X\times X$
with bi-transitivity without changing the class $\wm $, since the
only way that $(X\times X,T\times T)$ could be bi-transitive but
not transitive is if it has isolated points; but then $X$ has isolated
points and it is easy to check that in this case $(X\times X,T\times T)$
is not bi-transitive.

In particular, this implies that $\wm \subseteq \trans $, since $X$
is a factor of $X\times X$. The following is well known:

\begin{lem}
$\wm \subseteq \ttrans $. 
\end{lem}
\begin{proof}
Suppose that $(X,T)\in \wm $ but $X$ is not totally transitive.
For some $k$ there is a proper open subset $U\subseteq X$ such that
$T^{-k}U\subseteq U$ and $U$ is not dense in $X$. We may assume
that $U\cap T^{-i}U=\emptyset $ for $i\geq 1$; for if $U\cap T^{-1}U\neq \emptyset $
replace $U$ by $U\cap T^{-1}U$, and if $U\cap T^{-2}U\neq \emptyset $
replace $U$ by $U\cap T^{-2}U$, and so on; after $k-1$ steps we
will have $U$ as desired. Let $V=\cup _{i=0}^{k-1}T^{-i}(U\times U)$.
Then $T^{-1}V\subseteq V$. On the other hand $V$ is not dense in
$X\times X$ because it does not intersect $U\times T^{-1}U$.
\end{proof}
One consequence of this is that, since $\all _{\ttrans }$ isn't dense
in $\all _{\trans }$, neither is $\all _{\wm }$.

\begin{thm}
$\wm $ is a $G_{\delta }$ subset of $\all $ and is dense in $\all _{\ttrans }$. 
\end{thm}
\begin{proof}
The proof that $\wm $ is a $G_{\delta }$ is similar to the proof
that $\trans $ is a $G_{\delta }$ (theorem \ref{thm:trans-is-Gdelta-in-S}).
Density again follows from theorem \ref{thm:dense-isomorphism-class}
and the existence of minimal weakly mixing Cantor systems.
\end{proof}

\subsection{Strong Mixing}

A dynamical system $(Y,T)$ is \emph{strongly mixing} if, for every
two open nonempty sets $U,V\subseteq Y$, there are only finitely
many integers $n$ for which $U\cap T^{-n}V=\emptyset $. The class
of strong mixing systems is denoted by $\sm $

\begin{lem}
A transitive system $(Y,T)$ is strongly mixing if and only if for
every nonempty open set $U\subseteq Y$ there are only finitely many
$n$ for which $U\cap T^{-n}U=\emptyset $.
\end{lem}
\begin{proof}
One direction is trivial. For the other, suppose $U,V\subseteq Y$
are nonempty open sets. By transitivity we have $W=U\cap T^{-k}V\neq \emptyset $
for some $k\in \mathbb{N}$. By assumption, $W\cap T^{-n}W\neq \emptyset $
for all but finitely many $n\in \mathbb{Z}$. But $W\cap T^{-n}W\neq \emptyset $
implies $U\cap T^{-n+k}V\neq \emptyset $.
\end{proof}
In the group of automorphisms of a Lebesgue space the measure-theoretically
weak mixing systems are generic while the measure-theoretically strong
mixing systems are exotic. By analogy one would expect that the strong
mixing systems are exotic in $\all _{\ttrans }$. Surprisingly they
are not:

\begin{thm}
\label{thm:strong-mixing-is-generic}$\sm $ is generic in $\all _{\ttrans }$. 
\end{thm}
We do not know whether $\all _{\sm }$ is a $G_{\delta }$. 

In order to prove this we will use the following approximation of
strong mixing. Let $\varepsilon >0$. We say that a system $(X,T)$
is \emph{$\varepsilon $-strongly-mixing} if there exists an integer
$N$ such that, for every $k>N$ and every $x\in X$, it holds that
$B_{\varepsilon }(x)\cap T^{-k}(B_{\varepsilon }(x))\neq \emptyset $. 

Clearly $(X,T)$ is strong mixing if and only if it is $1/n$-strong-mixing
for every $n\in \mathbb{N}$. Let \[
V_{n}=\{X\in \all \, :\, X\textrm{ is }\frac{1}{n}-\textrm{strongly mixing}\}\]
 In order to prove theorem \ref{thm:strong-mixing-is-generic} it
suffices to show that the each $V_{n}$ contains an open dense subset
of $\all _{\ttrans }$. This follows from

\begin{prop}
\label{prop:strong-mixing-approximation}Let $\varepsilon >0$ and
$X\in \all _{\ttrans }$. For every $\delta >0$ there is a system
$Y\in \all _{\ttrans }$ with $d(X,Y)<\delta $ and an $r>0$ such
that any totally transitive system $Z$ with $d(Y,Z)<r$ is $\varepsilon $-strongly-mixing. 
\end{prop}
We will break the proof into a sequence of lemmas. First, we establish
a symbolic condition for $\varepsilon $-strong-mixing:

\begin{lem}
\label{lem:epsilon-strong-mixing}A system $X\in \all $ is $\varepsilon $-strong-mixing
if and only if there is a finite set $L\subseteq Q^{*}$of words appearing
in $X$ and an integer $N$ such that 
\end{lem}
\begin{enumerate}
\item Every subword of $X$ of length $1+2/\varepsilon $ is an $\varepsilon $-subword
of a member of $L$. 
\item For every $a\in L$ and every $k>N$ there is some $x\in X$ in which
$a$ $\varepsilon $-appears at index $0$ and at index $k$.
\end{enumerate}
\begin{proof}
For $x\in X$, note that $B_{\varepsilon }(x)\cap \sigma ^{-k}(B_{\varepsilon }(x))\neq \emptyset $
if and only if there is some $x'\in X$ with $d(x,x')<\varepsilon $
and $d(x,\sigma ^{k}x')<\varepsilon $. This in turn is equivalent
to the fact that the subwords of $x'$ of length $1+2/\varepsilon $
appearing at indices $0$ and $k$ differ by at most $\varepsilon $
from the central subword of $x$ of the same length. The lemma now
follows by compactness. 
\end{proof}
The main step in the proof of proposition \ref{prop:strong-mixing-approximation}
is the construction of $Y$. Fix $X\in \all _{\ttrans }$ and $\varepsilon ,\delta >0$.
We will construct $Y$ which is $\delta $-close to $X$ by symbolic
approximation as described in lemma \ref{prop:symbolic-approximation}.
Our aim is to perform the construction in such a way that $Y$ satisfies
the hypothesis of lemma \ref{lem:epsilon-strong-mixing}, and furthermore
every system close enough to $Y$ does as well.

Fix $\lambda >0$ to be determined later. Select two $\lambda $-pseudo-periods
in $Y$, say $\overline{a}=a_{1},\ldots ,a_{N}$ and $\overline{b}=b_{1},\ldots ,b_{M}$,
with $a_{1}=b_{1}$ and lengths $M,N$ relatively prime, such that

\renewcommand{\theenumi}{\alph{enumi}}

\begin{enumerate}
\item The collection \[
\{a_{k}\, :\, \frac{1}{\varepsilon }<k<N-\frac{1}{\varepsilon }\}\]
 is $\lambda $-dense in $Y$. 
\item For any sub-sequence $\overline{u}$ of $\overline{a}$ and any $0\leq k<N$
there is an index $i$ with $i=k(\mod N)$ such that $\overline{u}$
occurs in $\overline{b}$ at $i$. 
\end{enumerate}
\renewcommand{\theenumi}{\arabic{enumi}}

To construct $\overline{a}$, first select a transitive point $a_{1}$
and set $a_{k}=\sigma ^{k}a_{1}$, then choose $N$ so that (a) holds.
Next, enumerate the possibilities for $\overline{u},k$ as in (b)
and realize them in $\overline{b}$ one at a time, using the total
transitivity of $Y$ to ensure that any possible alignment can be
achieved given enough time. Extend $\overline{b}$ as needed to ensure
its length is prime relative to $N$.

Now $\overline{b}$ is much longer than $\overline{a}$. Choose $K$
so that the $K$-time concatenation $\overline{a}^{K}=\overline{a}\, \overline{a}\ldots \overline{a}$
is much longer than $\overline{b}$, say $KN>10M$, and such that
$KN$ and $M$ are relatively prime. Write $\overline{c}=\overline{a}^{K}$;
note that $\overline{c}$ is also a $\lambda $-dense $\lambda $-pseudo-period
starting at $a_{1}$. We assume that the language $\{\overline{b},\overline{c}\}$
has the unique parsing property; see definition \ref{def:unique-parsing}
and the remarks following it.

Set $a=\theta (\overline{a})\, ,\, b=\theta (\overline{b})\, ,\, c=\theta (\overline{c})$
(so $c$ is the $K$-times concatenation of $a$). By \ref{prop:symbolic-approximation},
any symbolic system constructed from $\{b,c\}$ will be $2\sqrt{\lambda }$-close
to $X$. Let $Y$ be such a totally transitive system, satisfying
the additional constraint that $b$ not appear in it twice consecutively.
Such a system exists since $KN,M$ are relatively prime.

\begin{lem}
If $\lambda $ was chosen small enough (in a manner independent of
$Y$) then $Y$ is $\varepsilon $-strongly-mixing. 
\end{lem}
\begin{proof}
Let $L$ be the set of subwords of $a$ of length $1+2/\varepsilon $.
We will claim that the two conditions in lemma \ref{lem:epsilon-strong-mixing}
hold with respect to $L$. 

Condition (1) follows from (a) assuming $\lambda $ was chosen small
enough; we omit the details.

We turn to (2). We claim that for every $w\in L$ and $n>KN$ there
is a $y\in Y$ such that $w$ appears in $y$ at indices $0$ and
$n$.

For let $n>KM$ and let $y\in Y$ be any transitive point for $y$.
Since $Y$ is totally transitive both $b$ and $c$ appear in $y$
(otherwise $y$, and hence $Y$, would be periodic). For convenience,
assume $b$ appears at index $0$. Consider the $M$-block in $y$
at location $n$. If this block is made up entirely of concatenations
of $a$'s then we are done: $w$ appears with period $N$ in $y$
between $n$ and $n+M$, and if its alignment is $k\mod N$ then we
can find the occurrence of $w$ in $b$ with this alignment, say at
index $j$ in $b$, and then $j,j+n$ are occurrences of $w$ in $y$
we were looking for.

Similarly, if the $M$-block starting in $y$ at index $-n$ is made
up entirely of concatenations of $a$ we are done.

If the argument above fails to produce the pair of $w$'s we wanted,
then it must be because there is an occurrence of $b$ at position
$n+n'$ in $y$ for some $-M\leq n'\leq M$ (or at $-n+n'$; this
case is dealt with similarly). Let $u_{1}u_{2}\ldots $ and $v_{1}v_{2}\ldots $
be the unique parsing of $y$ into $\{b,c\}$ words starting at index
$0$ and $n+n'$ respectively, so $u_{j},v_{j}\in \{b,c\}$ for all
$j$. By assumption $u_{1}=v_{1}=b$. If $u_{j}=v_{j}$ for all $j$
then the point $y$ would be periodic, which is impossible since $y$
generates $Y$ and $Y$ is totally transitive. So let $j$ be the
first index such that $u_{j}\neq v_{j}$. Then $u_{j-1}=v_{j-1}$
and they cannot be equal to $b$, because if both were equal to $b$
then the next block in both cases would necessarily be $c$ (recall
that $b$'s do not occur consecutively in $y$) and we get $u_{j}=v_{j}$,
a contradiction. So we have either $u_{j-1}u_{j}=c\, c=a^{K}a^{K}$
and $v_{j-1}v_{j}=cb=a^{K}b$, or else $u_{j-1}u_{j}=a^{K}b$ and
$v_{j-1}v_{j}=a^{K}a^{K}$. In either case since the beginning of
$u_{j},v_{j}$ are exactly $n+n'$ apart, and since the length of
$a^{K}$ is several times that of $b$ we are back in the situation
from the previous paragraph and can find two occurrences of $w$ separated
by $n$, as desired.
\end{proof}
To complete the proof of the proposition, we claim that if $Z$ is
totally transitive and is close enough to $Y$ then $Z$ enjoys the
same kind of properties as $Y$. To be precise,

\begin{lem}
Let $Y\in \all _{\ttrans }$ be constructed from $\{b,c\}$. For $\lambda $
small enough (in a manner independent of $Y$), if $Z$ is totally
transitive and sufficiently close to $Y$, then $Z$ is $10\varepsilon $-strongly-mixing. 
\end{lem}
\begin{proof}
For $Z$ close enough to $Y$ the projection $\tau _{Y}$ is defined
on $Z$. Furthermore since $\{b,c\}$ has the unique parsing property
when $Z$ is close enough to $Y$, we have that for any $z\in Z$
its projection $\tau _{Y}(z)$ is a concatenation of $b,c$'s with
no two $b$'s appearing consecutively. Also, since $Y$ is nontrivial
when $Z$ is close enough to $Y$, we see that $\tau _{Y}(Z)$ is
nontrivial. Therefore $\tau _{Y}(Z)$ is a totally transitive system
(since it is a factor of the totally transitive system $Z$) constructed
from $\{b,c\}$ and satisfying the constraint that the word $bb$
never appear, and hence by the previous lemma, $\tau _{Y}(Z)$ is
$\varepsilon $-strongly-mixing.

Now use the fact that if $Z$ is close enough to $Y$ then $d_{\infty }(\tau _{Y}(z),z)<\varepsilon $
for every $z\in Z$. Using lemma \ref{lem:epsilon-strong-mixing}
one sees that this implies that $Z$ is $10\varepsilon $-strongly
mixing for all $Z$ sufficiently close to $Y$. 
\end{proof}
This completes the proof of proposition \ref{prop:strong-mixing-approximation}
and that strong mixing is generic in $\all _{\ttrans }$. Theorem
\ref{thm:transitive-systems}(\ref{thm:totally-transitive-systems-strong-mixing})
now follows from the correspondence theorem.

\section{\label{sec:prime-and-doubly-minimal-systems}Minimal Self Joinings}

In this section by the \emph{orbit} of a point $x$ we will mean the
full orbit of $x$, i.e. the set $\{T^{k}x\, :\, k\in \mathbb{Z}\}$.
For a system $(X,T)$ we will denote the action of $T\times T$ on
$X\times X$ simply by $T$, so $T(x',x'')=(Tx',Tx'')$.

A topological system $(Y,S)$ has \emph{minimal self-joinings} if
it is infinite and for every $(x',x'')\in X\times X$, either $x',x''$
are on the same orbit (i.e. there is a $n\in \mathbb{Z}$ with $T^{n}x'=x''$)
or else the orbit of $(x',x'')$ is dense in $X\times X$. Such systems
are also said to be \emph{doubly minimal.}

There is a standard one-one correspondence between factors of $(X,T)$
and closed invariant equivalence relations of $(X\times X,T)$: A
factor $Y$ of $X$ given by a map $\varphi :X\rightarrow Y$ corresponds
to the closed equivalence relation $\{(x',x'')\in X\times X\, :\, \varphi (x')=\varphi (x'')\}$.

When $X$ has minimal self-joinings the only subsystems of $X\times X$
are the trivial subsystem $X\times X$ and the graphs \[
D_{n}=\{(x,T^{n}x)\, :\, x\in X\}\]
 for $n\in \mathbb{N}$. Now, $D_{n}$ is not an equivalence for $n\neq 0$:
otherwise for every $x\in X$, the fact that $(x,T^{n}x)\in D_{n}$
would imply (by symmetry of the equivalence relation $D_{n}$) that
$(T^{n}x,x)\in D_{n}$ and thus by definition of $D_{n}$ we get $T^{2n}x=x$.
Since this holds for every $x\in X$ we deduce that every point of
$X$ has period $2n$, which is impossible because $X$ has minimal
self-joinings (proof: since $X$ is infinite, we can choose two periodic
points not on the same orbit. This pair violates the definition of
double minimality).

Consequently, when $X$ has minimal self-joinings, the only closed
equivalence relations of $X\times X$ are the entire space $X\times X$,
which corresponds to the factor map taking all of $X$ to a single
fixed point, or the diagonal $D_{0}=\{(x,x)\, :\, x\in X\}$, which
corresponds to the factor given by the identity map on $X$. Thus
the only factors of a system with minimal self-joinings are the trivial
factors. Systems with this property are called \emph{prime}; the simplest
example being periodic systems of prime period. 

The notion of double minimality comes from a similar notion in the
measure preserving category, where examples were first constructed
by Rudolph \cite{R-msj}, and shown to be exotic by del Junco in \cite{dJ}.
However it was recently shown by Ageev that primeness in the measure-theoretic
category is generic \cite{A00}.

Prime topological systems were first constructed by Furstenberg, Keynes
and Shapiro \cite{FKS}. The first topological system with minimal
self-joinings was constructed by J. King in \cite{K}; Later B. Weiss
showed in \cite{W-dbmin} that every ergodic system $(Y,\mathcal{B},\mu ,T)$
with zero entropy has a realization as an invariant measure on a topological
system with minimal self-joinings. This implies that in some sense
topological double minimality is a common phenomenon. Double minimality
is a common phenomenon in our setup as well:

\begin{thm}
\label{thm:Double-minimality-is-generic}Double minimality is generic
in $\all _{\ttrans }$
\end{thm}
We will use the following combinatorial fact which is a watered-down
version of lemma 2 from \cite{W-dbmin}:

\begin{lem}
\label{lem:approximately-random-sets}For any integer $L$ and every
large enough $N$ there is a set $I\subseteq \{1,\ldots ,N\}$ such
that, if $i\in I$, then $i+k\notin I$ for $1\leq k\leq L$ and for
every $N/10\leq k\leq 9N/10$ it holds that $(I+k)\cap I\neq \emptyset $.
Such a set $I$ is called \emph{approximately random}.
\end{lem}
We also use the following simple observation:

\begin{lem}
\label{lem:distinct-orbits}Let $b=b(1)\ldots b(M)\in Q^{*}$ and
$c=c(1)\ldots c(N)\in Q^{*}$ and assume $2M/5<N<M/2$. Let $x,y\in Q^{\mathbb{Z}}$
be non-periodic points constructed from $b,c$ such that $c$ doesn't
appear twice consecutively in $x,y$. Then either $x=\sigma ^{n}y$
for some $n\in \mathbb{Z}$, or else there are occurrences $u,v$
of $b$ in $x,y$ whose alignment differs by some $k$ in the range
$(2M/10,8M/10)$.
\end{lem}
\begin{proof}
Since $x,y$ are constructed from $L$ we may parse $x,y$ as \begin{eqnarray*}
x & = & \ldots u_{-1}u_{0}u_{1}u_{2}\ldots \\
y & = & \ldots v_{-1}v_{0}v_{1}v_{2}\ldots 
\end{eqnarray*}
where $u_{i},v_{i}\in \{b,c\}$ (note that the parsing may not be
unique, so we choose one such parsing). 

Suppose there is no pair of $b$'s in $x,y$ whose alignment differs
by an integer in the range $(2M/10,8M/10)$. We will show that $x=\sigma ^{n}y$
for some $n$.

We claim that if $u_{i}=c$ for some $i$ then there is a unique $j$
such that $v_{j}=c$ and $u_{i},v_{j}$ overlap. Indeed, there can
obviously be at most one such $j$, since occurrences of $c$ are
separated by $b$'s and $b$ is longer than $c$. Now if the statement
were false then $u_{i}$ overlaps with one or two occurrences of $b$
in $y$. By the restriction that $cc$ does not appear in $x$, we
see that $u_{i-1}u_{i}u_{i+1}=bcb$, and one can check directly that
no matter what configuration of $b's$ occurs in $y$ opposite $u_{i}$,
there is a $b$ in $y$ opposite one of the words $u_{i-1}$ or $u_{i+1}$
at a shift in the range $(4M/5,M/2)$, contrary to our assumption.

Similarly, if $v_{i}=c$ then there is a unique occurrence of $c$
in the parsing of $x$ overlapping $v_{i}$.

Now suppose $u_{i}u_{i+1}\ldots u_{i+k}=cbb\ldots bbc$. The distance
between the end of $u_{i}$ and the beginning of $u_{i+k}$ is a multiple
of $M$. From the above it is clear that there are occurrences $v_{j},v_{j+m}$
of $c$ in $y$ opposite $u_{i},u_{i+k}$ respectively. There can
be no occurrence of $c$ in $y$ between $v_{j}$ and $v_{j+m}$ because
this would imply an occurrence in $x$ between $u_{i}$ and $u_{i+k}$,
contrary to our assumption. Thus $kM-2N\leq mM\leq kM+2N$. Since
$N<M/2$ we must have $k=m$. In particular the offset of $v_{j}$
from $u_{i}$ is equal to the offset of $u_{i+k}$ from $v_{j+k}$
and $u_{i}u_{i+1}\ldots u_{i+k}=v_{j}v_{j+1}\ldots v_{j+k}$. 

Similarly, suppose $u_{i}=c$ and $u_{k}=b$ for all $k>i$. The same
argument shows that there is a $j$ such that $v_{j}=c$ overlaps
$u_{i}$ and $v_{k}=b$ for all $k>j$. A similar statement is true
if $u_{k}=b$ for all $k<i$.

If $u_{i}=v_{i}=b$ for every $i$ then clearly $x,y$ are shifts
of each other. Otherwise pick an occurrence of $c$ in $x$, which
we may assume is $u_{0}$. Find the occurrence of $c$ in $y$ opposite
it, which we assume if $v_{0}$. Now look to the right; repeated application
of the above shows that the right-infinite subword of $x$ starting
at $u_{0}$ equals the right-infinite subword of $y$ starting at
$v_{0}$. A similar statement holds for left-infinite subwords, so
$y$ is a shift of $x$.
\end{proof}
The main step in the proof of theorem \ref{thm:Double-minimality-is-generic}
is the following construction.

\begin{lem}
Let $X\in \all _{\ttrans }$. For every $\varepsilon ,\delta $ there
exists a totally transitive symbolic system $Y=Y(X,\varepsilon ,\delta )$
with $d(X,Y)<\varepsilon $ and such that for every $y',y''\in Y$
either $y'=\sigma ^{n}y''$ for some $n\in \mathbb{Z}$ or else the
orbit of $\{\sigma ^{n}(y',y'')\, :\, n\in \mathbb{Z}\}$ is $\delta $-dense
in in $Y\times Y$. Furthermore, if $Z$ is close enough to $Y$ then
$\tau _{Y}$ takes $Z$ into $Y$. 
\end{lem}
\begin{proof}
Let $X,\varepsilon ,\delta $ be given. Choose a small positive $\eta $
and construct an $\eta $-dense $\eta $-pseudo-period $\overline{a}=a_{1},\ldots ,a_{R}$
in $X$. Now form a very long $\eta $-pseudo orbit $\overline{b}=b_{1},\ldots ,b_{M}$
in $X$ for which the set of indices $I\subseteq \{1,\ldots ,M\}$
at which $\overline{a}$ occurs has the property that $I\cap (I+k)\neq \emptyset $
for every $M/10\leq k\leq 9M/10$ (such a $\overline{b}$ exists by
total transitivity and lemma \ref{lem:approximately-random-sets}).
Let $1\leq i,j\leq R$ and $\frac{1}{5}M\leq k\leq \frac{4}{5}M$;
then there are occurrences of $a_{i}$ and $a_{j}$ in $\overline{b}$
whose alignment differs by $k$. To see this note that assuming $M>10R$
we can find, by approximate randomness, two occurrences of $\overline{a}$
in $\overline{b}$ whose alignment differs by $k+i-j$. 

The significance of the last paragraph is that if we take two copies
of $\overline{b}$ and shift one of them right or left by a distance
between $2M/10$ and $8M/10$, then every pair $a_{i},a_{j}$ occurs
in the pair $\overline{b}$ and $\sigma ^{k}\overline{b}$. Since
$\overline{a}$ is $\eta $-dense, the collection of pairs $(a_{i},a_{j})$
occurring in $(\overline{b},\sigma ^{k}\overline{b})$ is $\eta $-dense
in $X\times X$.

Choose another $\eta $-pseudo-orbit $c=c_{1},\ldots ,c_{N}$ with
$c_{1}=b_{1}$ and such that $4M/10<N<M/2$ and $M,N$ relatively
prime (we can do this because we may assume $M$ as large as necessary).
Write $b,c$ for the finite words $b=\theta (\overline{b)}$ and $c=\theta (\overline{c})$.
Let $Y\in \all $ be the system constructed from the language $\{b,bcb\}$
containing all infinite concatenations of $b,c$ in which the word
$cc$ does not appear. This is a strongly mixing system since $K,N$
are relatively prime (it's essentially a mixing subshift of finite
type). We have by lemma \ref{prop:symbolic-approximation} that $d(X,Y)<2\sqrt{\eta }$,
which can be made $<\varepsilon $.

We may assume that $\{b,c\}$ has the unique parsing property (definition
\ref{def:unique-parsing}) so $\{b,bcb\}$ do as well; this implies
that $\tau _{Y}(Z)\subseteq Y$ for $Z$ close enough to $Y$ (lemma
\ref{lem:symbolic-retract} and what follows).

It remains to show that every $(y',y'')\in Y\times Y$ not on the
same orbit have a $\delta $-dense orbit under $Y\times Y$. Applying
lemma \ref{lem:distinct-orbits} we see if $y',y''$ are not on the
same orbit then there are occurrences of $b$ in $y',y''$ whose alignment
difference is in the range $(2M/10,8M/10)$. Now $y'=\theta (\overline{y}')$
and $y''=\theta (\overline{y}'')$ for $\eta $-pseudo-orbits $\overline{y'},\overline{y}''$
made by concatenating $\overline{b}$'s and $\overline{c}$'s; and
therefore by our assumptions, every pair $(a_{i},a_{j})$ occurs in
the pair $\overline{y}',\overline{y}''$. But since $\{(a_{i},a_{j})\}_{1\leq i,j\leq R}$
is $\eta $-dense in $X\times X$ this means that the orbit of $(y',y'')=\theta (\overline{y}',\overline{y}'')$
is $\delta $-dense in $Y\times Y$, assuming $\eta $ small enough;
as required.
\end{proof}
{}

\begin{proof}
(of theorem \ref{thm:Double-minimality-is-generic}) As in section
\ref{sec:symbolic-approximation} let \[
\rho _{Y}=\min \{d(c',c'')\, :\, c',c''\in \pi _{0}(Y)\, ,\, c'\neq c''\}\]
 and recall that if $Y$ is symbolic then $\rho _{Y}>0$, and if $d(Z,Y)<\rho _{Y}$
then $\tau _{Y}$ is defined on $Z$ (lemma \ref{lem:symbolic-retract}).
Also write\[
\widetilde{\rho }_{Y}(Z)=\max _{z\in Z}\min _{y\in Y}d(z,y)\]
 For every $Z$ close enough to $Y$ we have $\widetilde{\rho }_{Y}(Z)<\frac{1}{2}\rho _{Y}$,
and then $\tau _{Y}$ is defined on $Z$, and furthermore $d_{\infty }(z,\tau _{Y}(z))<\frac{1}{2}\widetilde{\rho }_{Y}(Z)$
for all $z\in Z$. 

For $X\in \all _{\ttrans }$ and $\varepsilon ,\delta >0$ let $Y=Y(X,\varepsilon ,\delta )$
be as in the lemma. Fix $n\in \mathbb{N}$ and choose $r=r(X,\varepsilon ,\delta )<\delta $
in such a way that \[
Z\in B_{r}(Y)\; \Rightarrow \; \widetilde{\rho }_{Y}(Z)<\min \{\frac{1}{4}\rho _{Y},\frac{1}{n}\}\]
 this guarantees that $\tau _{Y}$ is defined on any $Z\in B_{r}(Y)$.
We may assume also that the projection of $\tau _{Y}(Z)$ of $Z$
is nontrivial, since we may assume $X$, and hence $Y$, are nontrivial. 

Write $B(X,\varepsilon ,\delta )=B_{r(X,\varepsilon ,\delta )}(Y(X,\varepsilon ,\delta ))$.
Set \[
U_{n}=\bigcup _{X\in \ttrans }\bigcup _{\varepsilon >0}B(X,\varepsilon ,\frac{1}{n})\]
 and \[
D=\bigcap _{n\in \mathbb{N}}U_{n}\]
 $D$ is a dense $G_{\delta }$ subset of $\all _{\ttrans }$ (each
$U_{n}$ is open and it is dense in $\all _{\ttrans }$ since since
$Y(X,\varepsilon ,\delta )\in \all _{\ttrans }$ for $X\in \all _{\ttrans }$).
We will show that if $Z\in D\cap \all _{\ttrans }\cap \all _{\mn }$
then $Z$ has minimal self-joinings; since $\mn $ is a dense $G_{\delta }$
in $\all _{\ttrans }$ this completes the proof.

Let $Z\in D\cap \all _{\ttrans }\cap \all _{\mn }$ and for $n\in \mathbb{N}$
let $Y_{n}=(X_{n},\varepsilon _{n},1/n)\, ,\, r_{n}=r(X_{N},\varepsilon _{n},1/n)$
such that $Z\in \cap B_{r_{n}}(Y_{n})$. Fix $n$ and $z',z''\in Z$
and let $y'=\tau _{Y_{n}}(z')$ and $y''=\tau _{Y_{n}}(z'')$. Distinguish
two cases: 
\begin{enumerate}
\item If $y',y''$ are not on the same orbit in $Y_{n}$ then the orbit
of $(y',y'')$ is $1/n$ dense in $Y_{n}\times Y_{n}$. For an arbitrary
pair $(w',w'')\in Z\times Z$ there is some $m$ such that \[
d(\sigma ^{m}(y'),\tau _{Y}(w'))<\frac{1}{n}\; \textrm{and}\; d(\sigma ^{m}(y''),\tau _{Y}(w''))<\frac{1}{n}\]
 Now since \[
d_{\infty }(y',z')<\widetilde{\rho }_{Y_{n}}(Z)\; \textrm{and}\; d_{\infty }(y'',z'')<\widetilde{\rho }{}_{Y_{n}}(Z)\]
 and since $\widetilde{\rho }{}_{Y_{n}}(Z)<1/n$ we can combine these
inequalities and get \[
d(\sigma ^{m}(y'),\sigma ^{m}(z'))<\frac{2}{n}\; \textrm{and}\; d(\sigma ^{m}(y''),\sigma ^{m}(z''))<\frac{2}{n}\]
 Similarly we have\[
d(w',\tau _{Y}(w'))<\frac{1}{n}\; \textrm{and}\; d(w'',\tau _{Y}(w''))<\frac{1}{n}\]
 from which we conclude\[
d(\sigma ^{m}(z'),w')<\frac{3}{n}\; \textrm{and}\; d(\sigma ^{m}(z''),w'')<\frac{3}{n}\]
 Since $w',w''$ were arbitrary the orbit of $(z',z'')$ in $Z\times Z$
is $3/n$ dense in $Z\times Z$. 
\item Otherwise, suppose that $y',y''$ are on the same orbit, so $y'=\sigma ^{k(n)}y''$
for some integer $k(n)$. Then since \[
d_{\infty }(y',z')\leq \widetilde{\rho }_{Y_{n}}(Z)<\min \{\frac{1}{2}\rho _{Y},\frac{1}{n}\}\]
 and the same is true of $y'',z''$ we see that $z'$ is $\min \{\frac{1}{2}\rho _{Y},\frac{1}{n}\}$-uniformly
close to $\sigma ^{k(n)}z''$. 
\end{enumerate}
If (1) holds for infinitely many $n$ then the orbit of $(z',z'')$
is dense in $Z\times Z$.

Otherwise (2) holds for infinitely many $n$. If there is some $m\in \mathbb{N}$
such that $k(n)=m$ for infinitely many $n$, we see that $z',\sigma ^{m}(z'')$
are $1/n$-uniformly close for infinitely many $n$, and so $z'=\sigma ^{m}(z'')$,
i.e. $z',z''$ are on the same orbit.

In the alternative case $k(n)$ would take on infinitely many values,
and in particular at least to values. Suppose for simplicity that
$k(1)\neq k(2)$ and that $\widetilde{\rho }_{Y_{1}}(Z)\leq \widetilde{\rho }_{Y_{2}}(Z)$.
We have than\[
d_{\infty }(z',\sigma ^{k(1)}(z''))<\frac{1}{2}\widetilde{\rho }_{Y_{1}}(Z)\; ,\; d_{\infty }(z',\sigma ^{k(2)}(z''))<\frac{1}{2}\widetilde{\rho }_{Y_{2}}(Z)\]
 and using the assumption $\widetilde{\rho }_{Y_{1}}(Z)\leq \widetilde{\rho }_{Y_{2}}(Z)$
we get\[
d_{\infty }(z'',\sigma ^{k(2)-k(1)}(z''))<\widetilde{\rho }_{Y_{2}}(Z)\]
 Therefore \[
\tau _{Y_{2}}(z'')=\sigma ^{k(2)-k(1)}(\tau _{Y_{2}}(z''))\]
 Since $Z$ is minimal the point $z''$ is a transitive point for
$Z$, so $\tau _{Y}(z'')$ is transitive for $Y$, and we conclude
that $Y$ is periodic with period $k(2)-k(1)$. But $Z$ is totally
transitive, so cannot have periodic factors; a contradiction. 
\end{proof}
This completes the proof the double minimality is generic in $\all _{\ttrans }$.
Theorem \ref{thm:transitive-systems}(\ref{thm:totally-transitive-systems-msj})
now follows from the correspondence theorem.

\section{\label{sec:open-questions}Closing Comments and Open Questions}

We briefly mention some problems and extensions which arise in connection
with this work.

Our original motivation for this work was our interest in $\all $
as a universal space of topological systems. That it is equivalent
to the homeomorphism group of $\cantor $ was a surprise and has proved
useful in the study of the latter. However there are many interesting
subsets of $\all $ which cannot be attacked by our methods.

A prime candidate is the space of connected systems in $\all $, which
is universal for connected systems. This is a closed space in $\all $,
hence Polish. What is generic there? More specifically, what is the
relation between the transitive, totally transitive and rigid connected
systems? What is generic in each? Is there a zero-one law? 

Another interesting question is, what is the connection between the
space of connected systems in $\all $ and the homeomorphism groups
of {}``nice'' manifolds? 

Throughout this paper we have been working with subsystems of the
shift space $Q^{\mathbb{Z}}$ where $Q$ is the Hilbert cube. It may
be verified that except for in the proof that $\all $ is universal,
the only property we have relied on has been that $Q$ has no isolated
points. For any perfect compact metric space $\Delta $ our proofs
work for the space $\all (\Delta )$ of closed shift-invariant subsystems
of $\Delta ^{\mathbb{Z}}$.

This is no longer true when $\Delta $ has isolated points. Most interesting
is the symbolic case, e.g. $\Delta =\{0,1\}$, so $\Delta ^{\mathbb{Z}}=\{0,1\}^{\mathbb{Z}}$
is the full shift on two symbols. Theorem \ref{sub:genericity-of-odometer}
is no longer true; indeed the only symbolic odometers are finite cycles.
Instead we have that the finite cycles form a dense $G_{\delta }$
(this also shows that the correspondence theorem is false here). Another
theorem which requires modification is theorem \ref{thm:dense-isomorphism-class};
to obtain an analogue of it, one must take entropy into account. Otherwise,
though, theorem \ref{thm:totally-transitive-systems} remains the
same in this space.

A problem we have only partly settled is the relation between genericity
in the space of measures and genericity of topological realizations
of them. We have provided a partial answer to this in theorem \ref{thm:trans-category-correspondence}.
It is possible that the hypothesis that the property be a $G_{\delta }$
is too strong. See the remarks after the proof of theorem \ref{thm:trans-category-correspondence}.

Finally, very little is known about genericity for actions of other
groups. Here part of our work carries over: the correspondence theorem
and its proof remain valid for such actions. More precisely consider
a countable group $\Gamma $ and the shift space $Q^{\Gamma }$ along
with the shift action of $\Gamma $ on $Q^{\Gamma }$ given by $(gx)(h)=x(hg)$
for $g,h\in \Gamma $. The correspondence theorem then relates genericity
in the space $\all _{\Gamma }$ of shift-invariant subsystems of $Q^{\Gamma }$
to genericity in the space of representations by homeomorphisms of
$\Gamma $ on the cantor set, with an appropriate topology, and similarly
for suitable subspaces.

However, theorems \ref{thm:transitive-systems},\ref{thm:totally-transitive-systems}
do not appear to have good analogues in this more general setting.
The following example is instructive. Consider the case of the group
$\Gamma =\mathbb{Z}$, and let $X\subseteq \{a,b\}^{\mathbb{Z}}$
be a shift of finite type, $a,b\in Q$. If $Y$ is close enough to
$X$ in $\all _{\mathbb{Z}}$, then by lemma \ref{lem:retract-embedding}
$Y$ has a nontrivial factor in $\{a,b\}^{\mathbb{Z}}$ in which each
sufficiently long word belongs to $X$; this implies there is a non-trivial
factor of $Y$ embedded in $X$. If $Y$ is prime, for instance, then
$Y$ embeds into $X$. In conjunction with our other results about
the genericity of prime systems and the density of the isomorphism
class of minimal totally transitive systems, this means that the subsystems
of mixing SFTs are a very large family. The last observation is of
course not new. The point is that when one goes to $\Gamma =\mathbb{Z}^{2}$
this argument cannot be true because there are minimal totally transitive
$\mathbb{Z}^{2}$ shifts of finite type $X$ (Moses \cite{M89}).
By the same reasoning as above, such a system has a neighborhood in
$\all _{\mathbb{Z}^{2}}$ consisting entirely of extensions of this
system, so the isomorphism class of any system not extending $X$
does not have a dense isomorphism class. Contrast this with theorem
\ref{thm:dense-isomorphism-class}.

\bibliographystyle{plain}
\addcontentsline{toc}{section}{\refname}\bibliography{bib}

\begin{thebibliography}{10}

\bibitem{A00}
O.~N. Ageev.
\newblock The generic automorphism of a {L}ebesgue space conjugate to a
  {$G$}-extension for any finite abelian group {$G$}.
\newblock {\em Dokl. Akad. Nauk}, 374(4):439--442, 2000.

\bibitem{AGW06}
E.~Akin, Eli Glasner, and B.~Weiss.
\newblock Generically there is but one homeomorphism of the cantor set.
\newblock {\em preprint, http://www.arxiv.org/abs/math.DS/0603538}, 2006.

\bibitem{AHK03}
Ethan Akin, Mike Hurley, and Judy~A. Kennedy.
\newblock Dynamics of topologically generic homeomorphisms.
\newblock {\em Mem. Amer. Math. Soc.}, 164(783):viii+130, 2003.

\bibitem{A78}
Steve Alpern.
\newblock Generic properties of measure preserving homeomorphisms.
\newblock In {\em Ergodic theory (Proc. Conf., Math. Forschungsinst.,
  Oberwolfach, 1978)}, volume 729 of {\em Lecture Notes in Math.}, pages
  16--27. Springer, Berlin, 1979.

\bibitem{AP02}
Steve Alpern and V.~S. Prasad.
\newblock Properties generic for {L}ebesgue space automorphisms are generic for
  measure-preserving manifold homeomorphisms.
\newblock {\em Ergodic Theory Dynam. Systems}, 22(6):1587--1620, 2002.

\bibitem{BDK}
Sergey Bezugly, Anthony~H. Dooley, and Jan Kwiatkowski.
\newblock Topologies on the group of homeomorphisms of a cantor set.
\newblock 2004.

\bibitem{CP}
J.~R. Choksi and V.~S. Prasad.
\newblock Approximation and {B}aire category theorems in ergodic theory.
\newblock In {\em Measure theory and its applications (Sherbrooke, Que.,
  1982)}, volume 1033 of {\em Lecture Notes in Math.}, pages 94--113. Springer,
  Berlin, 1983.

\bibitem{dJ}
Andr{\'e}s del Junco.
\newblock Disjointness of measure-preserving transformations, minimal
  self-joinings and category.
\newblock In {\em Ergodic theory and dynamical systems, I (College Park, Md.,
  1979--80)}, volume~10 of {\em Progr. Math.}, pages 81--89. Birkh\"auser
  Boston, Mass., 1981.

\bibitem{F}
Harry Furstenberg.
\newblock Disjointness in ergodic theory, minimal sets, and a problem in
  {D}iophantine approximation.
\newblock {\em Math. Systems Theory}, 1:1--49, 1967.

\bibitem{FKS}
Harry Furstenberg, Harvey Keynes, and Leonard Shapiro.
\newblock Prime flows in topological dynamics.
\newblock {\em Israel J. Math.}, 14:26--38, 1973.

\bibitem{GK-01law}
Eli Glasner and Jonathan~L. King.
\newblock A zero-one law for dynamical properties.
\newblock In {\em Topological dynamics and applications (Minneapolis, MN,
  1995)}, volume 215 of {\em Contemp. Math.}, pages 231--242. Amer. Math. Soc.,
  Providence, RI, 1998.

\bibitem{GW-TRproperty}
Eli Glasner and Benjamin Weiss.
\newblock The topological {R}ohlin property and topological entropy.
\newblock {\em Amer. J. Math.}, 123(6):1055--1070, 2001.

\bibitem{H44}
Paul~R. Halmos.
\newblock Approximation theories for measure preserving transformations.
\newblock {\em Trans. Amer. Math. Soc.}, 55:1--18, 1944.

\bibitem{Hal-mixing}
Paul~R. Halmos.
\newblock In general a measure preserving transformation is mixing.
\newblock {\em Ann. of Math. (2)}, 45:786--792, 1944.

\bibitem{Hal-ergodicth}
Paul~R. Halmos.
\newblock {\em Lectures on ergodic theory}.
\newblock Chelsea Publishing Co., New York, 1960.

\bibitem{KR04}
Alexander~S. Kechris and Christian Rosendal.
\newblock Turbulence, amalgamation and generic automorphisms of homogeneous
  structures.
\newblock {\em preprint, http://www.arxiv.org/abs/math.LO/0409567}, 2004.

\bibitem{K}
Jonathan~L. King.
\newblock A map with topological minimal self-joinings in the sense of del
  {J}unco.
\newblock {\em Ergodic Theory Dynam. Systems}, 10(4):745--761, 1990.

\bibitem{M89}
Shahar Mozes.
\newblock Tilings, substitution systems and dynamical systems generated by
  them.
\newblock {\em J. Analyse Math.}, 53:139--186, 1989.

\bibitem{OU41}
J.~C. Oxtoby and S.~M. Ulam.
\newblock Measure-preserving homeomorphisms and metrical transitivity.
\newblock {\em Ann. of Math. (2)}, 42:874--920, 1941.

\bibitem{Ox}
John~C. Oxtoby.
\newblock {\em Measure and category}, volume~2 of {\em Graduate Texts in
  Mathematics}.
\newblock Springer-Verlag, New York, second edition, 1980.
\newblock A survey of the analogies between topological and measure spaces.

\bibitem{Rk-mixing}
V.~Rohlin.
\newblock A ``general'' measure-preserving transformation is not mixing.
\newblock {\em Doklady Akad. Nauk SSSR (N.S.)}, 60:349--351, 1948.

\bibitem{R}
D.~Rudolph.
\newblock Residuality and orbit equivalence.
\newblock In {\em Topological dynamics and applications (Minneapolis, MN,
  1995)}, volume 215 of {\em Contemp. Math.}, pages 243--254. Amer. Math. Soc.,
  Providence, RI, 1998.

\bibitem{R-msj}
Daniel~J. Rudolph.
\newblock An example of a measure preserving map with minimal self-joinings,
  and applications.
\newblock {\em J. Analyse Math.}, 35:97--122, 1979.

\bibitem{Wal81}
Peter Walters.
\newblock {\em An introduction to ergodic theory}, volume~79 of {\em Graduate
  Texts in Mathematics}.
\newblock Springer, 1981.

\bibitem{W-dbmin}
B.~Weiss.
\newblock Multiple recurrence and doubly minimal systems.
\newblock In {\em Topological dynamics and applications (Minneapolis, MN,
  1995)}, volume 215 of {\em Contemp. Math.}, pages 189--196. Amer. Math. Soc.,
  Providence, RI, 1998.

\end{thebibliography}

\end{document}